\documentclass{amsart}

\input xy 
\xyoption{all} 
\xyoption{v2} 
\xyoption{2cell} 

\newcommand{\Hom}{\mathbf{Hom}}

\renewcommand{\a}{\alpha}
\renewcommand{\b}{\beta}
\renewcommand{\c}{\gamma}
\renewcommand{\d}{\delta}

\newcommand{\equivalencerelation}{\sim}
\newcommand{\Aut}{\text{Aut}}
\newcommand{\Reals}{\mathbb R}
\renewcommand{\P}{{\mathcal P}}

\newcommand{\U}{\mathcal U} 
\newcommand{\Z}{\mathbb Z}
\newcommand{\C}{\mathbb C}

\newcommand{\cstar}{\C^\times}
\newcommand{\BGrb}{\mathbf{BGrb}} 
\newcommand{\XpiM}{X\stackrel{\pi}{\to}M} 

\theoremstyle{plain}
\newtheorem{theorem}{Theorem}[section]
\newtheorem{lemma}[theorem]{Lemma}
\newtheorem{proposition}[theorem]{Proposition}

\theoremstyle{definition}
\newtheorem{definition}[theorem]{Definition}

\theoremstyle{remark}
\newtheorem{example}{Example}[section]

\begin{document} 

\title{Bundle 2-Gerbes} 

\author{Daniel Stevenson} 

\thanks{The author acknowledges the support of the 
Australian Research Council} 
\address 
{Department of Pure Mathematics \\ 
University of Adelaide        \\ 
Adelaide, SA 5005             \\ 
Australia} 
\email{dstevens@maths.adelaide.edu.au} 


\begin{abstract} 
We make the category $\BGrb_M$ of bundle gerbes on a manifold $M$  
into a $2$-category by providing $2$-cells in the form 
of transformations of bundle gerbe morphisms.   
This description of $\BGrb_M$ as a $2$-category is used to 
define the notion of a bundle $2$-gerbe.  To every 
bundle $2$-gerbe on $M$ is associated a class in 
$H^4(M;\Z)$.  We define the notion of a bundle $2$-gerbe 
connection and show how this leads to a closed, 
integral, differential $4$-form on $M$ which represents 
the image in real cohomology of the class in 
$H^4(M;\Z)$.  Some examples of bundle $2$-gerbes 
are discussed, including the bundle $2$-gerbe 
associated to a principal $G$ bundle $P\to M$.  It 
is shown that the class in $H^4(M;\Z)$ associated to this 
bundle $2$-gerbe coincides with the first Pontryagin 
class of $P$ --- this example was previously considered 
from the point of view of $2$-gerbes 
by Brylinski and McLaughlin.  
\end{abstract} 

\UseTwocells 
\UseHalfTwocells 

\maketitle 

\section{Introduction} 
\label{sec:one} 

Recently there has been interest in developing 
higher dimensional analogues of line bundles --- 
so-called $p$-gerbes or $p$-line bundles --- which 
realise classes in $H^{p+1}(M;\Z)$ for a manifold $M$.  Part 
of the motivation for this comes from physicists, 
who wish to interpret closed $p$-forms with integral 
periods on $M$ as a generalised curvature 
of a bundle-like object on $M$.  A first step towards 
this goal was taken in the book \cite{Bry} of 
Brylinski, who developed a theory of differential 
geometry for \emph{gerbes}.  Gerbes were orginally 
introduced (in a very general setting) 
by Giraud in \cite{Gir} for the purposes 
of developing a degree $2$ non-abelian 
cohomology theory.  The theory described by Brylinski 
allows one to realise classes in $H^3(M;\Z)$ 
as equivalence classes of (abelian) gerbes.  
Murray in \cite{Mur} invented the notion 
of a \emph{bundle gerbe}.  Bundle gerbes are simpler 
objects than gerbes but still provide a 
geometric realisation of $H^3(M;\Z)$.  The 
theory of gerbes and bundle gerbes has proved 
to be very useful tool:  in  
\cite{CarMur} and \cite{CarMicMur} the authors studied 
anomalies in quantum field theory with the aid of bundle 
gerbes, Hitchin in \cite{Hit} has used the theory 
of gerbes in his study of mirror symmetry, while 
Brylinski has made extensive applications of gerbes --- 
one example is his use of gerbes in \cite{Bry1} to give 
an interpretation of Beilinson's regulator maps in 
algebraic $K$-theory.  
        
In \cite{BryMac} and \cite{BryMac1} the authors 
constructed a canonical $2$-gerbe associated to a
principal $G$ bundle $P\to M$ where $G$ is a 
compact, simple, simply connected Lie group.  $2$-gerbes, 
introduced by Breen in \cite{Bre}, are higher dimensional 
analogues of gerbes.  Breen 
used $2$-gerbes to study three dimensional non-abelian 
sheaf cohomology, however there is a certain class 
of $2$-gerbes --- $2$-gerbes \emph{bound} by the sheaf 
of abelian groups $\underline{\C}^{\times}_M$ --- 
that give rise to classes in $H^4(M;\Z)$ via 
the exponential isomorphism $H^3(M;\underline{\C}^{\times}_M) 
= H^4(M;\Z)$.  This is the class of $2$-gerbes studied 
by Brylinski and McLaughlin.  They show that the 
canonical $2$-gerbe associated to the principal bundle $P$ 
has class in $H^4(M;\Z)$ equal to $p_1$, the first 
Pontryagin class of $P$.  

We shall consider here a related geometric object, 
the bundle $2$-gerbe.  Bundle $2$-gerbes were originally 
introduced in \cite{CarMurWan} --- we shall use a 
modification of the definition used there.  A bundle $2$-gerbe 
is a quadruple of manifolds $(Q,Y,X,M)$ where 
$(Q,Y,X^{[2]})$ is a bundle gerbe \cite{Mur} over the fibre 
product $X^{[2]}$.  We also require that there is a 
bundle $2$-gerbe product.  In fact this requires two 
product structures, the first of which is a product on $Y$, which 
on the fibres takes the form 
$Y_{(x_2,x_3)}\times Y_{(x_1,x_2)}\to Y_{(x_1,x_3)}$  
for points $x_1,x_2$ and $x_3$ all lying in the same fibre.  
There is also a product in $Q$ covering this product 
on $Y$, and which commutes with the bundle gerbe 
product on $(Q,Y,X^{[2]})$.  This product on $Q$ 
satisfies a certain associativity condition.  One 
can associate to a bundle $2$-gerbe $(Q,Y,X,M)$ a $\cstar$ 
valued \v{C}ech $3$-cocycle $g_{ijkl}$ representing a 
class in $H^4(M;\Z)$.  One can also develop the notion 
of a bundle $2$-gerbe connection and a $2$-curving for a 
bundle $2$-gerbe connection in an analogous manner to 
\cite{Mur} and show that a bundle $2$-gerbe equipped 
with such structures has a $4$-curvature.  This is 
a closed, integral differential $4$-form on $M$ which is a 
representative in $H^4(M;\Reals)$ for the image, 
in real cohomology, of the class in $H^4(M;\Z)$ 
defined by the cocycle $g_{ijkl}$.  

There is a naturally arising bundle $2$-gerbe $Q$ 
associated to a principal $G$ bundle $P$ on 
$M$ where $G$ is as above.  If one calculates 
the \v{C}ech cocycle $g_{ijkl}$ associated to $Q$ 
then one recovers the results of \cite{BryMac} and 
\cite{BryMac1} giving an explicit cocycle formula 
for the first Pontryagin class of $P$.  

In outline then this paper is as follows.  In 
Section~\ref{sec:two} we review the theory of 
bundle gerbes from \cite{Mur}.  In Section~\ref{sec:three} 
we discuss a gluing or `descent' construction 
for line bundles from \cite{Bry}.  In Section~\ref{sec:four} 
we explain how to make the category of bundle 
gerbes on a manifold $M$ into a $2$-category 
by adding $2$-cells in the form of transformations 
of bundle gerbe morphisms.  This allows us 
in Section~\ref{sec:five} to `categorify' 
the definition of a bundle 
gerbe, so as to define a bundle $2$-gerbe.  The 
relationship of bundle $2$-gerbes with 
bicategories \cite{Ben} is also examined 
here.  This is also preparation for 
Section~\ref{sec:six} where an example of a 
bundle $2$-gerbe --- the tautological bundle 
$2$-gerbe --- is introduced via the homotopy  
bigroupoid of a space.  A \v{C}ech $3$-class 
is associated to a bundle $2$-gerbe in 
Section~\ref{sec:seven} and a de Rham 
representative for this class is defined 
in Section~\ref{sec:eight} via the notion 
of a bundle $2$-gerbe connection.  In 
Section~\ref{sec:nine} the example of a 
bundle $2$-gerbe associated to a principal 
$G$-bundle is discussed and, using the work 
of Brylinski and McLaughlin, it is shown that 
the $4$-class of this bundle $2$-gerbe 
coincides with the first Pontryagin class 
of the bundle.  In Sections~\ref{sec:ten} 
and ~\ref{sec:eleven} we discuss higher descent 
properties of bundle $2$-gerbes and define the 
notion of a trivial bundle $2$-gerbe.  We 
finally show that a bundle $2$-gerbe is trivial 
if and only if its $4$-class vanishes.  
We will not discuss the relationship of bundle $2$-gerbes 
with $2$-gerbes, this will be done elsewhere \cite{Ste1}.  For some 
preliminary results in this 
direction one can consult \cite{Ste}.  

This work is clearly influenced by the ideas 
presented in \cite{BryMac} and \cite{BryMac1}. 
I am very grateful to Michael Murray for his 
supervision of my PhD thesis and for his 
help in the preparation of this paper.     

\section{Review of Bundle Gerbes}
\label{sec:two} 

Let $\pi\colon X\to M$ be a surjection admitting local 
sections.  Let $X^{[2]} = X\times_{M}X$ denote the 
fiber product of $X$ with itself over $M$ and let 
$X^{[p]} = X\times_{M}X\times_{M}\cdots \times_{M}X$ 
denote the $p$-fold such fiber product.   
We can form a simplicial manifold 
$X_{\bullet} = \{X_{p}\}$ with $X_{p} = X^{[p+1]}$ and 
the face and degeneracy operators $d_{i}$ and $s_{i}$ 
given by omitting the $i^{\text{th}}$ factor and 
repeating the $i^{\text{th}}$ factor respectively.  
Thus the face operators $d_{i}\colon X^{[p+1]}\to 
X^{[p]}$ are given by $d_i = \pi_i$ where  
$$
\pi_{i}(x_{1},\ldots,x_{p+1}) = (x_{1},\ldots,x_{i-1},
x_{i+1},\ldots,x_{p+1}) 
$$
for $i=1,\ldots,p+1$ and for $p=1,2,\ldots$.  
Recall from \cite{Mur} that a \emph{bundle gerbe} consists 
of a triple $(P,X,M)$ where $\pi\colon X\to M$ is a 
surjection admitting local sections and $P$ is a 
principal $\cstar$ bundle on $X^{[2]}$ with a product.  
This means that there is a $\cstar$ bundle isomorphism 
$$
m_{P}\colon \pi_{1}^{-1}P\otimes \pi_{3}^{-1}P \to 
\pi_{2}^{-1}P 
$$
covering the identity on $X^{[3]}$.  Here $\pi_{1}^{-1}
P\otimes \pi_{3}^{-1}P$ denotes the \emph{contracted 
product} of the $\cstar$ bundles $\pi_{1}^{-1}P$ and 
$\pi_{3}^{-1}P$ --- see \cite{Bry}.  Fiberwise the 
bundle gerbe product $m_{P}$ is a map 
$$
m_{P}\colon P_{(x_{2},x_{3})}\otimes P_{(x_{1},x_{2})} 
\to P_{(x_{1},x_{3})} 
$$
for $(x_{1},x_{2},x_{3}) \in X^{[3]}$ and we usually 
write $u_{23}u_{12}$ for $m_P(u_{23}\otimes u_{12})$ 
when $u_{23}\in P_{(x_2,x_3)}$ and $u_{12}\in P_{(
x_1,x_2)}$.  The bundle gerbe 
product $m_{P}$ is required to be associative in the 
following sense: whenever $u_{34} \in P_{(x_{3},x_{4})}$, 
$u_{23} \in P_{(x_{2},x_{3})}$ and $u_{12} \in P_{(x_{1},
x_{2})}$ for $(x_{1},x_{2},x_{3},x_{4}) \in X^{[4]}$ 
we have $u_{34}(u_{23}u_{12}) = (u_{34}u_{23})u_{12}$.  
When $M$ is understood we will frequently write 
$(P,X)$ or even $P$ for $(P,X,M)$.    

Recall that a bundle gerbe also has an \emph{identity} 
section; this is a section $e$ of $P$ over the diagonal 
$\Delta(X) = \{(x,x)|x\in X\}\subset X^{[2]}$ which 
behaves as an identity with respect to the bundle gerbe 
product.  So if $u\in P_{(x_1,x_2)}$ then 
we have $ue(x_1) = u = e(x_2)u$.  A bundle gerbe also 
has an inverse map $P\to inv^{-1}P$ where 
$inv\colon X^{[2]}\to X^{[2]}$ is the map which 
switches an ordered pair $(x_1,x_2)$, so $inv
(x_1,x_2) = (x_2,x_1)$.  We denote the image of $u\in 
P_{(x_1,x_2)}$ under $P\to inv^{-1}P$ by $u^{-1}$ 
--- this has all the desired 
properties: $uu^{-1} = e(x_2)$, $(uv)^{-1} = v^{-1}u^{-1}$ 
and so on.  Note also that we can identify $inv^{-1}
P$ with $P^{*}$, the $\cstar$ bundle $P$ with the action 
of $\cstar$ changed to its inverse.  For more details 
we refer to \cite{Mur}.  

Various operations can be performed on bundle gerbes; for 
example there is the notion of the \emph{pullback} 
$(f^{-1}P,f^{-1}X,N)$ of a bundle gerbe $(P,X)$ on $M$  
by a map $f\colon N\to M$.  One can also form the 
\emph{product} $(P\otimes Q,X\times_M Y)$ of two bundle 
gerbes $(P,X)$ and $(Q,Y)$ on $M$.  Given a bundle gerbe 
$(P,X)$ we can also form its \emph{dual} $(P^*,X)$.  
We refer to \cite{Mur} 
for more details on these constructions.  

Suppose $Q\to X$ is a principal $\cstar$ 
bundle on $X$ and $\pi\colon X\to M$ is a 
local-section-admitting surjection.  Let 
$P$ be the $\cstar$ bundle on $X^{[2]}$ with 
fibre 
\begin{equation} 
\label{eq:trivial bundle gerbe} 
P_{(x,y)} = \Aut_{\cstar}(Q_x,Q_y) 
\end{equation} 
at $(x,y)\in X^{[2]}$.  $Q$ has an associative 
product via composition of isomorphisms.  A 
bundle gerbe isomorphic to a bundle gerbe of the 
form~(\ref{eq:trivial bundle gerbe}) via an 
isomorphism preserving the bundle gerbe products 
is said to be \emph{trivial}.   
The notation $\d(Q) = \pi_1^{-1}Q\otimes 
\pi_2^{-1}Q^*$ is frequently used to 
denote the bundle gerbe~(\ref{eq:trivial 
bundle gerbe}).   

In \cite{Mur} the notion of a \emph{bundle gerbe connection} 
on a bundle gerbe $(P,X)$ was introduced.  Before 
we recall this notion it is useful to note (see 
\cite{CarMur}) that we 
can reformulate the definition of a bundle gerbe in 
terms of line bundles and line bundle isomorphisms 
by replacing the principal $\cstar$ bundle $P$ with its 
associated line bundle $L$.  Then $L$ has an associative 
product $m_L\colon \pi_1^{-1}L\otimes \pi_3^{-1}L\to 
\pi_2^{-1}L$ described in the same manner above.  A bundle 
gerbe connection on $P$ then is a connection $\nabla_L$ 
on $L$ which is compatible with the bundle gerbe product 
$m_L$ in the sense that 
$$
\pi_1^{-1}\nabla_L + \pi_3^{-1}\nabla_L = m_L^{-1}\circ 
\pi_2^{-1}\nabla_L \circ m_L.   
$$
It is easy to see that the curvature $F_{\nabla_L}$ 
of a bundle gerbe connection $\nabla_L$ satisfies 
$\d(F_{\nabla_L}) = 0$.  Here $\d\colon \Omega^{p}(X^{[q]}) 
\to \Omega^{p}(X^{[q+1]})$ is the map formed by adding 
the pullback maps $\pi_i^*$ with an alternating sign: 
$\d = \sum (-1)^i \pi_i^*$.  Therefore $\d$ commutes 
with the exterior derivative $d$ and, since the $\pi_i$ are face 
maps for a simplicial manifold, it follows  
that $\d^2 = 0$.   
Hence we have a complex 
\begin{equation} 
\label{eq:exact complex of forms} 
\Omega^p(M) \stackrel{\pi^*}{\to} \Omega^p(X) \stackrel{\d}
{\to} \Omega^p (X^{[2]}) \stackrel{\d}{\to} \cdots 
\stackrel{\d}{\to} \Omega^p (X^{[q]}) \stackrel{\d}{\to} \cdots  
\end{equation} 
It is a fundamental result of \cite{Mur} 
that the complex~(\ref{eq:exact complex of forms})  
has no cohomology as long as $M$ supports partitions of unity.  
Hence we can solve 
the equation $F_{\nabla_L} = \d(f)$ for some two form 
$f$ on $X$.  Following \cite{Mur} we call a choice of 
this two form $f$ a \emph{curving} for the bundle gerbe 
connection $\nabla_L$.  From the equation $F_{\nabla_L} 
= \d(f)$ we obtain $\d(df) = 0$ and hence $df = 
\pi^*(\omega)$ for some necessarily closed three form $\omega$ 
on $M$.  One can show that $\omega$ has integral 
periods and hence is a representative of the image 
in $H^3(M;\Reals)$ of 
a class in $H^3(M;\Z)$.  We call the 
three form $\omega$ the \emph{$3$-curvature} of the 
bundle gerbe connection $\nabla_L$ and curving $f$.  

One can associate to any bundle gerbe $P$ on $M$ a 
$\underline{\C}^{\times}_M$-valued \v{C}ech $2$-cocycle 
$g_{ijk}$ as described in \cite{Mur}.  $g_{ijk}$ is a 
representative of a characteristic  
class $DD(P)$ in $H^3(M;\Z)$ --- the Dixmier-Douady 
class of the bundle gerbe $P$.  The 
$3$-curvature $\omega$ of a bundle gerbe connection 
on $P$ is a representative for the image, in real 
cohomology, of $DD(P)$.  The Dixmier-Douady class 
has the following properties.  

\begin{proposition}[\cite{Mur}] 
The Dixmier-Douady class $DD(P)$ of a bundle gerbe 
$P$ on $M$ satisfies 
\begin{enumerate} 
\item $DD(P\otimes Q) = DD(P) + DD(Q)$ for bundle gerbes 
$P$ and $Q$ on $M$.  
\item $DD(P^*) = - DD(P)$ where $P^*$ is the dual of the 
bundle gerbe $P$.  
\item $DD(f^{-1}P) = f^* DD(P)$ where $f^{-1}P$ denotes the 
pullback of the bundle gerbe $P$ on $M$ by a map $f\colon 
N\to M$.  
\end{enumerate} 
\end{proposition} 

Recall from \cite{Mur} that a \emph{bundle gerbe morphism} 
$f\colon P \to Q$ between bundle gerbes $P = (P,X)$  
and $Q = (Q,Y)$ is a triple of maps $f = (\hat{f},
f,\phi)$ where $f\colon X\to Y$ is a map commuting with 
the projections $\pi_{X}\colon X\to M$, $\pi_{Y}\colon 
Y \to M$ and covering $\phi\colon M\to M$, while $\hat{f}
\colon P\to Q$ is a $\cstar$ bundle morphism covering the 
induced map $f^{[2]}\colon X^{[2]}\to Y^{[2]}$.  We will 
only be interested in the case where $\phi = \text{id}_{M}$.    
One could define an isomorphism of bundle gerbes $P$ and 
$Q$ to be a morphism of bundle gerbes $(\hat{f},f,\phi)
\colon P\to Q$ in which each map was an isomorphism, however 
it is not true that isomorphism classes of bundle gerbes 
are in a bijective correspondence with $H^3(M;\Z)$.  
Instead, one can consider the weaker notion of 
\emph{stable isomorphism} \cite{MurSte} of bundle 
gerbes and show that there is a bijection between 
stable isomorphism classes of 
bundle gerbes and $H^3(M;\Z)$.  
 
\section{The Generalised Clutching Construction} 
\label{sec:three} 

Recall the following result 
from \cite{Bry}.  
\begin{lemma}[\cite{Bry}]  
\label{lemma:descent} 
Suppose $\pi\colon X\to M$ is a surjection admitting 
local sections and that $P$ is a $\cstar$ bundle on 
$X$ together with an isomorphism $\phi\colon \pi_2^{-1}
P\to \pi_1^{-1}P$ which satisfies the descent cocycle 
condition  
\begin{equation} 
\label{eq:descent cocycle} 
\pi_1^{-1}\phi\circ \pi_3^{-1}\phi = \pi_2^{-1}\phi 
\end{equation} 
over $X^{[3]}$.  Then the $\cstar$ bundle $P$ descends to 
$M$, ie there is a $\cstar$ bundle $Q = D(P)$ on $M$ 
plus an isomorphism $\psi\colon P\to \pi^{-1}Q$ which 
is compatible with $\phi$.  The converse is also true.  
\end{lemma} 

The $\cstar$ bundle isomorphism $\phi$ above is called 
a \emph{descent isomorphism}.  Note that fiberwise $\phi$ is 
a map $P_{x_1}\to P_{x_2}$ and the descent cocycle 
condition~(\ref{eq:descent cocycle}) is simply that the diagram  
$$
\xymatrixrowsep{8pt}
\xymatrixcolsep{5pt} 
\diagram  
P_{x_1}  \ar[rr] \ar[dr] & &   P_{x_2} \ar[dl]  \\ 
& P_{x_3} &                 \\ 
\enddiagram  
$$    
commutes.  We give an example of this kind 
of formalism below.  

\begin{example} 
\label{ex:difference of two trivialisations} 
Suppose $(P,X)$ is a bundle gerbe 
on $M$ and suppose that there are 
two trivialisations $T_1$ and $T_2$ 
of $P$ on $X$.  Thus there exist 
isomorphisms $P = \d(T_1)$ and $P = 
\d(T_2)$ commuting with the respective 
bundle gerbe products.  It is easy to 
see that there is a trivialisation of the 
bundle $\d(T_1\otimes T_2^*)$ over 
$X^{[2]}$.  This corresponds to an 
isomorphism $\phi\colon \pi_1^{-1}
(T_1\otimes T_2^*) \to \pi_2^{-1}(T_1 
\otimes T_2^*)$ covering the identity on 
$X^{[2]}$.  Since the isomorphisms 
$P = \d(T_1)$ and $P = \d(T_2)$ 
commute with the bundle gerbe products 
on the respective bundle gerbes, 
one can show that $\phi$ satisfies 
the descent cocycle condition.  Hence 
the bundle $T_1\otimes T_2^*$ descends 
to a bundle $D$ on $M$, ie there is an 
isomorphism $T_1 = T_2\otimes \pi^{-1}D$ 
of bundles on $X$, where $\pi\colon X\to M$ 
denotes the projection.  
\end{example} 

There is the following strengthening of 
the above lemma \cite{Bry}: 
there is an equivalence of categories $D\colon 
\textbf{Desc} (\XpiM ) \to \textbf{Bund}_M$ between the  
so called \emph{descent} category 
$\textbf{Desc} (\XpiM )$ and the category 
of principal $\cstar$ bundles $\textbf{Bund}_M$ on 
$M$.  Here $\textbf{Desc} (\XpiM )$ 
is the category whose 
objects are pairs $(P,\phi)$ where $\phi\colon \pi_2^{-1}P
\to \pi_1^{-1}P$ is a descent isomorphism as above and 
whose arrows $(P,\phi)\to (Q,\psi)$ are $\cstar$ 
bundle isomorphisms $f\colon P \to Q$ compatible with $\phi$ 
and $\psi$, so the following diagram commutes: 
$$
\xymatrix{ 
\pi_2^{-1}P \ar[r]^-{\pi_2^{-1}f} \ar[d]_-{\phi} & 
\pi_2^{-1}Q \ar[d]^-{\psi}                         \\ 
\pi_1^{-1}P \ar[r]^-{\pi_1^{-1}f} & \pi_1^{-1}Q.  } 
$$
It is clear that the operation $D$ which associates 
the $\cstar$ bundle $D(P)$ on $M$ to a bundle $P$ on 
$X$ with a descent isomorphism $\phi$ extends to an 
operation on maps --- if $f\colon (P,\phi)\to (Q,\psi)$ 
then there is an induced map $D(f)\colon D(P)\to D(Q)$ --- 
and this operation is functorial with respect to 
composition of maps.  

One other point to note is that if we make $\textbf{Desc}
(\XpiM )$ and $\textbf{Bund}_M$ into monoidal categories 
via the contracted product $\otimes$ of $\cstar$ bundles, then 
the equivalence of categories $D\colon \textbf{Desc}(\XpiM )
\to \textbf{Bund}_M$ commutes with $\otimes$ up to natural 
isomorphism.  More specifically, we define a functor 
$\otimes \colon \textbf{Desc}(\XpiM )\times \textbf{Desc}(\XpiM ) 
\to \textbf{Desc}(\XpiM )$ by a map on objects given by 
$\otimes ((P,\phi),(Q,\psi)) = (P\otimes Q,\phi\otimes \psi)$ 
and by a map on arrows given by $\otimes (f,g) = f\otimes g$.  
Then there is a natural isomorphism between the functors 
bounding the following diagram: 
$$
\xymatrix{ 
\textbf{Desc}(\XpiM )\times 
\textbf{Desc}(\XpiM ) \ar[d]_-{\otimes} 
\ar[r]^-{D\times D} & \textbf{Bund}_M 
\times \textbf{Bund}_M \ar @2{->}[dl]  
\ar[d]^-{\otimes}            \\ 
\textbf{Desc}(\XpiM ) \ar[r]^-{D} 
& \textbf{Bund}_M.             }  
$$
Note that such an isomorphism amounts to an isomorphism 
$D(P)\otimes D(Q) \to D(P\otimes Q)$ which is natural 
with respect to maps.     

\section{The 2-Category of Bundle Gerbes} 
\label{sec:four} 
 
Given bundle gerbes $P = (P,X)$ and $Q = (Q,Y)$ 
together with a pair of bundle gerbe morphisms 
$f,g\colon P\to Q$ with $f = (\hat
{f},f)$ and $g = (\hat{g},g)$ let $\hat{D}_{f,g}$ 
denote the $\cstar$ bundle $(f,g)^{-1}Q$ on $X$.  
Therefore $\hat{D}_{f,g}$ has fibre $Q_{(f(x),g(x))}$ 
at $x\in X$.  We will construct a descent isomorphism 
$\phi_{f,g}\colon \pi_2^{-1}\hat{D}_{f,g}\to \pi_1^{-1}
\hat{D}_{f,g}$ for $\hat{D}_{f,g}$.  Suppose $v \in 
(\pi_2^{-1}\hat{D}_{f,g})_{(x_1,x_2)} = (\hat{D}_{f,g})_
{x_1}$.  Thus $v \in Q_{(f(x_1),g(x_1))}$.  
Choose $u \in P_{(x_1,x_{2})}$ and put $\phi_{f,g}(v) = 
\hat{g}(u)(v\hat{f}(u^{-1}))$.   
Notice that this is independent of the choice of $u \in 
P_{(x_1,x_2)}$.  $\phi_{f,g}$ is a 
descent isomorphism --- ie it satisfies  
$$
\pi_1^{-1}\phi_{f,g}\circ \pi_3^{-1}\phi_{f,g} = 
\pi_2^{-1}\phi_{f,g} 
$$
over $X^{[3]}$.  This is a consequence of the 
associativity of the bundle gerbe products on $P$ and 
$Q$.   
We have the following Lemma.  
\begin{lemma}[\cite{Ste}] 
\label{lemma:bg transformations} 
\begin{enumerate} 

\item Suppose $(P,X)$ and $(Q,Y)$ are 
bundle gerbes on $M$ and that   
there exist bundle gerbe morphisms 
$f\colon P\to Q$ and $g\colon P\to Q$.  
Then the $\cstar$ bundle $\hat{D}_{f,g} = 
(f,g)^{-1}Q$ on $X$ descends to a $\cstar$ 
bundle $D_{f,g} = D(\hat{D}_{f,g})$ on $M$.  

\item Suppose that $P$, $Q$, $f$ and 
$g$ are as above and that there is a third 
bundle gerbe morphism $h \colon P\to Q$.   
Then there is an 
isomorphism 
$$
D_{g,h}\otimes D_{f,g}  
 \simeq D_{f,h} 
$$
of $\cstar$ bundles on $M$.  

\item Suppose that $P$ and $Q$ are as above but 
now we have bundle gerbe morphisms $f,g,
h,k  \colon 
P\to Q$.  Then the following diagram of $\cstar$  
bundle isomorphisms on $M$ commutes: 
$$
\xymatrix{ 
D_{h,k}\otimes D_{g,h} 
\otimes D_{f,g}  
\ar[d] \ar[r] & D_{h,k}\otimes 
D_{f,h} \ar[d]                 \\ 
D_{g,k}\otimes D_{f,g}  
\ar[r] & D_{f,k}    } 
$$
where the isomorphisms  
are those of (2) above.  

\end{enumerate} 
\end{lemma} 
 
(2) of this lemma is proved by noticing that the bundle gerbe 
product on $Q$ gives an isomorphism $(g,h)^{-1}Q\otimes 
(f,g)^{-1}Q \to (f,h)^{-1}Q$ of $\cstar$ bundles on $X$ 
which commutes with the descent isomorphisms for $(g,h)^{-1}Q
\otimes (f,g)^{-1}Q$ and $(f,h)^{-1}Q$ respectively.  Therefore 
there is an induced isomorphism $D((g,h)^{-1}Q\otimes 
(f,g)^{-1}Q)\to D((f,h)^{-1}Q)$.  (3) of the lemma is proved 
similarly, using the associativity of the bundle gerbe 
product on $Q$, the functorality of the operation $D$, and 
the fact that $D$ commutes with $\otimes$ up to natural 
isomorphism.  
  
This Lemma suggests the following Definition. 
\begin{definition}[\cite{Ste}] 
\label{def:bg transformations} 
Let $(P,X)$ and $(Q,Y)$ be bundle gerbes on $M$.  
A \emph{transformation} $\theta\colon f\Rightarrow 
g$ between two bundle gerbe morphisms $f, 
g\colon P\to Q$ 
is a section 
of the $\cstar$ bundle $D_{f,g} = 
D(\hat{D}_{f,g})$ on $M$.  
\end{definition} 

We would like to form a category $\Hom(P,Q)$ 
associated to bundle gerbes $(P,X)$ and $(Q,Y)$ 
with the bundle gerbe morphisms $P\to Q$ as objects.  
Therefore we would like to be able to compose transformations 
between bundle gerbe morphisms.  A way to do this is 
suggested by the previous lemma.  Given bundle gerbe 
morphisms $f,g,h\colon P\to Q$ 
together with transformations $\theta\colon f
\Rightarrow g$ and $\lambda\colon g\Rightarrow 
h$ then we have the induced section $\lambda \otimes 
\theta$ of $D_{g,h}\otimes D_{f,g}$.  We define the 
composed transformation $\lambda\theta\colon f
\Rightarrow h$ to be the image of this section 
$\lambda\otimes \theta$ under the isomorphism $D_{g,h}
\otimes D_{f,g}\to D_{f,h}$.  By the lemma above this 
operation of composition is associative.  We can define 
an identity transformation $1_f\colon f\Rightarrow 
f$ by noticing that the identity section of the 
bundle gerbe $Q$ pullsback to define a section $\hat{1}_f$ 
of $(f,f)^{-1}Q$ which is compatible with the descent 
isomorphism for $\hat{D}_{f,f} = (f,f)^{-1}Q$.  Therefore 
it descends to a section $1_f$ of $D_{f,f}$ and it is straightforward 
to check that this acts as an identity.  

The case where the manifold $M$ is a point illuminates the preceding 
discussion.  
In this case a bundle 
gerbe over a point becomes a $\cstar$ groupoid --- 
ie a groupoid such that the automorphism groups of 
each object of the groupoid are isomorphic to $\cstar$.  
Following \cite{Mur} we define the $\cstar$ groupoid $\textbf{Gr}(P)$  
associated to a bundle gerbe $(P,X,M)$ when the manifold $M$ 
is restricted to a point $m_0 \in M$ as follows.  We let the objects of 
the groupoid $\textbf{Gr}(P)$ be the points of $X_{m_0}$ where 
$X_{m_0} = \pi^{-1}(m_0)$.  Given two 
points of $X_{m_0}$, $x_1$ and $x_2$, we define the set of 
arrows $\text{Hom}(x_1 , x_2)$ from $x_1$ to $x_2$ in 
$\textbf{Gr}(P)$ to be the points of the fiber $P_{(x_1 , 
x_2)}$.  Composition of arrows in $\textbf{Gr}(P)$ is then 
provided by the bundle gerbe product on $P$ and 
the identity arrow from a point $x$ to itself is 
provided by the identity section $e(x)$ of $P$ 
evaluated at the point $x$.  Since inverses exist 
in $P$ every arrow is invertible and it is not hard to 
see that $\textbf{Gr}(P)$ is a $\cstar$ groupoid.  Thus we have 
a family of $\cstar$ groupoids, indexed by the points 
of $M$.  It is in this sense that a bundle gerbe is a 
`bundle of groupoids'.    

It is not hard to see that in this case, when $M$ is 
restricted to a point, a bundle gerbe morphism $f\colon P\to 
Q$ induces a functor $f\colon \textbf{Gr}(P)
\to \textbf{Gr}(Q)$ (the important point here is that 
$\hat{f}$ preserves the bundle gerbe products on 
$P$ and $Q$).  Suppose that we are given a second 
bundle gerbe morphism $g\colon P\to Q$ 
and a transformation $\theta\colon f\Rightarrow 
g$.  So $\theta$ is a section of the $\cstar$ 
bundle $D_{f,g} = D(\hat{D}_{f,g})$ on $M$ and 
hence lifts to a section $\hat{\theta}$ of the 
$\cstar$ bundle $\hat{D}_{f,g} = (f,g)^{-1}Q$ on 
$X$.  It follows from the definition of $\hat{D}_{f,g}$ 
that we have the following isomorphism of $\cstar$ 
bundles on $X^{[2]}$: 
$$
\psi\colon \pi_1^{-1}\hat{D}_{f,g}\otimes (f^{[2]})^{-1}Q 
\stackrel{\simeq}{\to} (g^{[2]})^{-1}Q\otimes \pi_2
^{-1}\hat{D}_{f,g}.  
$$
It also follows that the section $\hat{\theta}$ 
of $\hat{D}_{f,g}$ is compatible with this isomorphism 
in the sense that $\psi(\hat{\theta}(x_2)\otimes 
\tilde{f}(u)) = \tilde{g}(u)\otimes \hat{\theta}(x_1)$ 
where $u \in P_{(x_1 , x_2 )}$ and $\tilde{f}\colon 
P \to (f^{[2]})^{-1}Q$ and $\tilde{g}\colon P\to 
(g^{[2]})^{-1}Q$ are induced by $\hat{f}$ and 
$\hat{g}$ respectively.  When we restrict $M$ to a point 
$m_0 \in M$, this is exactly the condition that 
$\hat{\theta}$ defines a natural transformation 
(in fact a natural isomorphism) between the functors 
$f$ and $g$.  

We would like to define a \emph{$2$-category} $\BGrb_M$ 
whose objects are the bundle gerbes $P$ on $M$.  We refer 
to \cite{KelStr} for the definition of a $2$-category (see also 
Section~\ref{sec:five}).  We take as the objects of $\BGrb_M$ 
the bundle gerbes $P$ on $M$, and given two bundle gerbes 
$P$ and $Q$ on $M$, we define the category $\Hom(P,Q)$ 
as above.  Thus the objects of $\Hom(P,Q)$ ($1$-arrows 
of $\BGrb_M$) are the bundle gerbe morphisms $P\to Q$ 
and the arrows of $\Hom(P,Q)$ ($2$-arrows of $\BGrb_M$) 
are the transformations $\theta\colon f\Rightarrow g$.  
We need to define a composition functor 
$$
m\colon \Hom(Q,R)\times \Hom(P,Q) 
\to \Hom(P,R).  
$$
It is clear how to define the action of $m$ on 
$1$-arrows: if $g\colon Q\to R$ and $f
\colon P\to Q$ are bundle gerbe morphisms, then 
we put $m(g,f) = g\circ f$.  
It is not so clear how to define the action of 
$m$ on $2$-arrows.  However we have the following result 
from \cite{Ste}.  

\begin{lemma} 
Suppose we are given three bundle gerbes $(P,X)$, 
$(Q,Y)$ and $(R,Z)$ on $M$ together with bundle gerbe 
morphisms $f_1,f_2\colon P\to Q$ and 
$g_1,g_2\colon Q\to R$.   
Then we have the following isomorphism 
of $\cstar$ bundles on $M$:
$$
D_{g_1\circ f_1,g_2\circ f_2} \simeq D_{f_1,f_2}\otimes 
D_{g_1,g_2}.  
$$
\end{lemma} 

This Lemma suggests a way to define the action 
of $m$ on $2$-arrows.  Suppose $\theta\colon f_1
\Rightarrow f_2$ is a transformation between 
bundle gerbe morphisms $(P,X)\to (Q,Y)$ and that $\lambda
\colon g_1\Rightarrow g_2$ is a transformation 
of bundle gerbe morphisms $(Q,Y)\to (R,Z)$.  Then $\theta$ 
and $\lambda$ lift to sections $\hat{\theta}$ and 
$\hat{\lambda}$ of the $\cstar$ bundles $\hat{D}_{f_1,
f_2} = (f_1,f_2)^{-1}Q$ and $\hat{D}_{g_1,g_2} = 
(g_1,g_2)^{-1}R$ on $X$ and $Y$ respectively.  
$\hat{g}_2$ induces an isomorphism $\tilde{g}_2\colon 
Q\to (g_2^{[2]})^{-1}R$, so if $x\in X$ then 
$\tilde{g}_2(\hat{\theta}(x)) \in R_{(g_2\circ f_1(x), 
g_2\circ f_2(x))}$.  Let $f_1^{-1}\hat{\lambda}$ 
denote the section of the pullback bundle $f_1^{-1}
\hat{D}_{g_1,g_2} = (g_1\circ f_1,g_2\circ f_1)^{-1}R$ 
on $X$.  Then if $x\in X$, $f_1^{-1}\hat{\lambda}(x) 
\in R_{(g_1\circ f_1(x),g_2\circ f_1(x))}$ and so 
$\tilde{g}_2(\hat{\theta})f_1^{-1}\hat{\lambda}(x) 
\in R_{(g_1\circ f_1(x),g_2\circ f_2(x))}$.  It is easy 
to check that $\tilde{g}_2(\hat{\theta})f_1^{-1}\hat{
\lambda}$ commutes with the descent isomorphism for 
$\hat{D}_{g_1\circ f_1,g_2\circ f_2} = (g_1\circ f_1,
g_2\circ f_2)^{-1}R$ and therefore descends to a section 
$\theta\circ \lambda$ of $D_{g_1\circ f_1,g_2\circ f_2}$.  

Note that fiberwise, ie 
regarding bundle gerbes as being bundles of groupoids, 
this is simply the operation of composing natural transformations 
in the $2$-category $\textbf{Cat}$ with categories 
as objects, functors as $1$-arrows and natural transformations 
as $2$-arrows --- recall that if we have categories 
$\textbf{C}$, $\textbf{D}$ and $\textbf{E}$ together with 
functors $F_1,F_2\colon \textbf{C}\to \textbf{D}$ and 
$G_1,G_2\colon \textbf{D}\to \textbf{E}$ plus natural 
transformations $\a\colon F_1\Rightarrow F_2$ and 
$\b\colon G_1\Rightarrow G_2$ then one can define the 
composed natural transformation $\b\circ \a\colon G_1\circ 
F_1\Rightarrow G_2\circ F_2$.  

We need to show that the action of $m$ on $2$-arrows 
is functorial.  Suppose that we have bundle 
gerbes $P$, $Q$ and $R$, bundle gerbe morphisms 
$f_1,f_2,f_3\colon P\to Q$, 
$g_1,g_2,g_3\colon Q\to R$ and 
transformations between them as pictured in the 
following diagram 
$$
\xymatrix{ 
P \rruppertwocell^{f_1}{\ \ \theta_{12}} \ar[rr]_(0.35){f_2}  
\rrlowertwocell_{f_3}{\ \ \theta_{23}} & & Q 
\rruppertwocell^{g_1}{\ \ \lambda_{12}} \ar[rr]_(0.35){g_2}  
\rrlowertwocell_{g_3}{\ \ \lambda_{23}} & & R. } 
$$
To show that $m$ is a functor we need to show that the 
two different ways of composing $2$-arrows coincide --- 
ie $(\lambda_{23}\lambda_{12})\circ (\theta_{23}\theta_{12}) 
= (\lambda_{23}\circ \theta_{23})(\lambda_{12}\circ \theta_{12})$.  
It is sufficient to show that 
$$
\tilde{g}_3(\hat{\theta}_{23})f_2^{-1}\hat{\lambda}_{23}
\tilde{g}_2(\hat{\theta}_{12})f_1^{-1}\hat{\lambda}_{12} 
= \tilde{g}_3(\hat{\theta}_{23}\hat{\theta}_{12})
f_1^{-1}(\hat{\lambda}_{23}\hat{\lambda}_{12}).  
$$
Since $\hat{\lambda}_{23}$ is compatible with the 
descent isomorphisms for the bundle $\hat{D}_{g_2,g_3}$, 
we have $\hat{\lambda}_{23}(y_2)\hat{g}_2(u) = 
\hat{g}_3(u)\hat{\lambda}_{23}(y_1)$ for $(y_1,y_2) 
\in Y^{[2]}$ and $u\in Q_{(y_1,y_2)}$.  Therefore 
$f_2^{-1}\hat{\lambda}_{23}\tilde{g}_2(\hat{\theta}_{12}) 
= \tilde{g}_3(\hat{\theta}_{12})f_1^{-1}\hat{\lambda}_{23}$, 
which establishes the equation above.    
One can also check 
that the functor $m$ is associative and that identity 
$1$-arrows and identity $2$-arrows behave as they should 
with respect to composition by $m$.  Hence we have the 
following proposition.  

\begin{proposition}[\cite{Ste}] 
There is a $2$-category $\BGrb_M$ whose 
objects are bundle gerbes $P$ on $M$, $1$-arrows 
are bundle gerbe morphisms $P\to Q$ and whose 
$2$-arrows are transformations between bundle gerbe 
morphisms with the composition laws given as 
above.  
\end{proposition} 
           
\section{Simplicial Bundle Gerbes and Bundle 2-Gerbes} 
\label{sec:five} 

We use our description of $\BGrb_{M}$ as a $2$-category 
to define the notion of a \emph{simplicial bundle gerbe} 
on a simplicial manifold $X = \{X_p\}$.  We are motivated 
by Brylinski and McLaughlin's definitions of a 
\emph{simplicial line bundle} (\cite{BryMac} Definition 5.1) 
and a \emph{simplicial gerbe} (\cite{BryMac} page 617).  We 
record here the definition of a simplicial line bundle.  

\begin{definition}[\cite{BryMac}] 
A \emph{simplicial line bundle} on a simplicial manifold 
$X_{\bullet} = \{X_p\}$ consists of the following data: 
\begin{enumerate} 

\item a line bundle $L\to X_1$ 
 
\item a non-vanishing section $s$ of the line 
bundle $\d(L)$ on $X_2$ where 
$$
\d(L) = d_0^{-1}L\otimes d_1^{-1}L^{*}\otimes 
d_2^{-1}L,  
$$
where $d_i\colon X_p \to X_{p-1}$ denote the 
face operators of the simplicial manifold 
$X_{\bullet} = \{X_p\}$.  

\item $s$ induces a non-vanishing section $\d(s)$ 
of the line bundle $\d \d(L)$ on $X_3$ where $\d \d(L)$ 
is defined by 
$$
\d \d(L) = d_0^{-1}\d(L)\otimes d_1^{-1}\d(L)^{*}\otimes 
d_2^{-1}\d(L)\otimes d_3^{-1}\d(L)^{*} 
$$
and $\d(s) = d_0^{-1}s\otimes d_1^{-1}s^{*}\otimes 
d_2^{-1}s\otimes d_3^{-1}s^{*}$.  Notice that as a result 
of the simplicial identities satisfied by the face 
operators $d_i\colon X_p \to X_{p-1}$ the line bundle 
$\d \d(L)$ is canonically trivialised.  We demand that 
$\d(s)$ matches this canonical trivialisation.  

\end{enumerate} 
\end{definition} 
 
Note the following consequences of this definition.  
\begin{list} 
{(\roman{enumi})}{\usecounter{enumi}} 

\item The non-vanishing section $s$ of $\d(L)$ defines 
a line bundle isomorphism $d_0^{-1}L\otimes d_2^{-1}L 
\to d_1^{-1}L$ covering the identity on $X_2$.  The 
coherency condition on $s$ is equivalent to this line 
bundle isomorphism satisfying an `associativity' 
condition on $X_3$.  

\item In the special case where the simplicial manifold 
$X_{\bullet} = \{X_p\}$ is the simplicial manifold  
associated to a surjection $\pi\colon X\to M$ which 
locally admits sections, then a simplicial line bundle 
on $X_{\bullet}$ recovers the definition of a bundle gerbe.  

\item Another important special case is when $X_{\bullet} = \{X_p\}$ 
is the simplicial manifold $NG$ associated to the classifying 
space of a Lie group $G$ (see \cite{Dup}).  Then a simplicial 
line bundle on $NG$ is the same thing as a central extension 
of $G$ by $\cstar$ (\cite{BryMac}). 

\end{list} 
 
We use the notion of a simplicial line bundle to motivate 
our definition of a simplicial bundle gerbe.   
To avoid cluttered notation later on, it is convenient 
to restrict attention to the simplicial 
manifold $X_{\bullet}$ associated to a surjection $\pi\colon 
X\to M$ which admits local sections.  This will not affect 
our results at all; everything we say will be true for 
an arbitrary simplicial manifold, however it is 
easier to state this for $X_{\bullet}$.  

We start with a bundle gerbe $(Q,Y,X^{[2]})$ on 
$X^{[2]}$.  We suppose there is a bundle gerbe morphism 
$m\colon \pi_1^{-1}Q\otimes \pi_3^{-1}Q\to 
\pi_2^{-1}Q$.  It is convenient to introduce 
some new notation (analogous to that used in 
\cite{Mur}) to avoid large, complicated diagrams.    
Let us denote by $Y\circ Y\to X^{[3]}$ the local-section-admitting 
surjection whose fiber at a point $(x_1,x_2,x_3)\in 
X^{[3]}$ is $Y_{(x_2,x_3)}\times Y_{(x_1,x_2)}$.  
So $Y\circ Y = \pi_1^{-1}Y\times_{X^{[3]}}\pi_3^{-1}Y$.   
Another way of looking at this is that $Y\circ Y$ is 
the restriction of $Y\times Y$ to $(Y\times Y)|_{X^{[2]}
\circ X^{[2]}}$ where $X^{[2]}\circ X^{[2]} = \{((x,y),
(y,z))|\ (x,y), (y,z) \in X^{[2]}\} = X^{[3]}$.  A 
point of $Y\circ Y$ is of the form $(y_{23},y_{12})$ 
where $y_{23} \in Y_{(x_2,x_3)}$ and $y_{12}\in 
Y_{(x_1,x_2)}$ for some point $(x_1,x_2,x_3)\in X^{[3]}$.     
Similarly let $Q\circ Q$ denote the restriction 
of $Q\otimes Q\to Y^{[2]}\times Y^{[2]}$ to 
$(Y\circ Y)^{[2]}\subset Y^{[2]}\times Y^{[2]}$.  
Thus $Q\circ Q =\pi_1^{-1}Q\otimes \pi_3^{-1}Q$ 
and has fiber $(Q\circ Q)_{((y_{23},y_{12}),(y_{23}',
y_{12}'))}$ at a point $((y_{23},y_{12}),
(y_{23}',y_{12}'))$ of $(Y\circ Y)^{[2]}$ equal 
to $Q_{(y_{23},y_{23}')}\otimes Q_{(y_{12},y_{12}')}$.  

By construction the triple $(Q\circ Q,Y\circ Y,X^{[3]})$ 
is a bundle gerbe --- the bundle gerbe $(\pi_1^{-1}Q
\otimes \pi_3^{-1}Q,\pi_1^{-1}Y\times_{X^{[3]}}\pi_3^{-1}Y,
X^{[3]})$.  The bundle gerbe morphism $m\colon 
\pi_1^{-1}Q\otimes \pi_3^{-1}Q\to \pi_2^{-1}Q$ is 
then a bundle gerbe morphism (also denoted $m$) 
$Q\circ Q\to Q$ covering the map $\pi_2\colon X^{[3]}
\to X^{[2]}$ sending 
a point $(x_1,x_2,x_3)$ of $X^{[3]}$ to the point 
$(x_1,x_3)$ of $X^{[2]}$.  Over $X^{[4]}$ we can define 
another bundle gerbe $(Q\circ Q\circ Q,Y\circ Y\circ Y,
X^{[4]})$ where $Y\circ Y\circ Y\to X^{[4]}$ is the 
local-section-admitting surjection with fiber 
$$
Y_{(x_3,x_4)}\times Y_{(x_2,x_3)}\times Y_{(x_1,x_2)} 
$$
over a point $(x_1,x_2,x_3,x_4)$.  
$Q^{\circ^3} = Q\circ Q\circ Q$ 
is defined in an analogous fashion to $Q\circ Q$ 
above.  The bundle gerbe morphism $m$ gives rise to 
two bundle gerbe morphisms $m_1, m_2 \colon 
Q^{\circ^3}\to Q$ which cover the map $X^{[4]}\to 
X^{[2]}$ which sends $(x_1,x_2,x_3,x_4)$ to $(x_1,x_4)$.  
We have $m_1 = (\hat{m}_1,m_1)$, $m_2 = 
(\hat{m}_2,m_2)$ where $m_1,m_2 \colon 
Y^{\circ ^3} = Y\circ Y\circ Y \to Y$ 
are given by $m_1(y_{34},y_{23},y_{12}) = m(m(y_{34},y_{23}),
y_{12})$ and    
$m_2(y_{34},y_{23},y_{12}) = m(y_{34},m(y_{23},
y_{12}))$, and $\hat{m}_1,\hat{m}_2\colon Q^{\circ^3}
\to Q$ are given by $\hat{m}_1(u_{34}\otimes u_{23}\otimes 
u_{12}) = \hat{m}(\hat{m}(u_{34}\otimes u_{23})\otimes 
u_{12})$ and $\hat{m}_2(u_{34}\otimes u_{23}\otimes u_{12}) 
= \hat{m}(u_{34}\otimes \hat{m}(u_{23}\otimes u_{12}))$, 
for $u_{ij} \in Q_{(y_{ij},y_{ij}^{'})}$.  We demand 
that there is a transformation of bundle gerbe morphisms 
$a\colon m_1\Rightarrow m_2$.  Recall that 
this means there is a section $\hat{a}$ of the $\cstar$ 
bundle $(m_1,m_2)^{-1}Q$ on $Y^{\circ^3}$ which 
descends to a section $a$ of the $\cstar$ bundle $A = 
D((m_1,m_2)^{-1}Q) = D_{m_1,m_2}$ on $X^{[4]}$.    

\begin{figure}  
$$ 
\xymatrix@!C=30pt{ 
& Q^{\circ^4}  \ar[rr]^-{m\circ 1\circ 1} 
\ar[dl]_-{1\circ 1\circ m} \ar[dr]^-{1\circ m\circ 1} 
& & Q^{\circ^3}  \ar[dr]^-{m\circ 1} 
\ar @2{->}[dl]_-{a\circ 1} &                                   \\ 
Q^{\circ^3}  \ar[dr]_-{1\circ m} &    & 
Q^{\circ^3}  \ar @2{->}[ll]_-{1\circ a} \ar[dl]_-{1\circ m} 
\ar[rr]^-{m\circ 1} &  & Q^{\circ^2}  \ar[dl]^-{m} 
\ar @2{->}[dlll]^-{a}                                             \\ 
& Q^{\circ^2}  \ar[rr]_-{\ \ \ m} & \ar @2{-}[d] & Q     &   \\ 
& Q^{\circ^4}  \ar[dl]_-{1\circ 1\circ m} 
\ar[rr]^-{\ \ \ \ \ \ \ \ \ \ m\circ 1\circ 1} & & Q^{\circ^3}  
\ar[dl]^-{1\circ m} \ar[dr]^-{m\circ 1} 
\ar @2{-}[dlll] &                                                 \\ 
Q^{\circ^3}  \ar[dr]_-{1\circ m} 
\ar[rr]_-{m\circ 1} & & Q^{\circ^2}  \ar[dr]_-{m} 
\ar @2{->}[dl]_-{a} & & Q^{\circ^2}  \ar @2{->}[ll]^-{a} 
\ar[dl]^-{m}                                               \\ 
& Q^{\circ^2}  \ar[rr]_-{m} & & Q                               } 
$$ 
\caption{Coherency condition for associator transformation} 
\label{fig:associativity coherence for bundle gerbes} 
\end{figure}
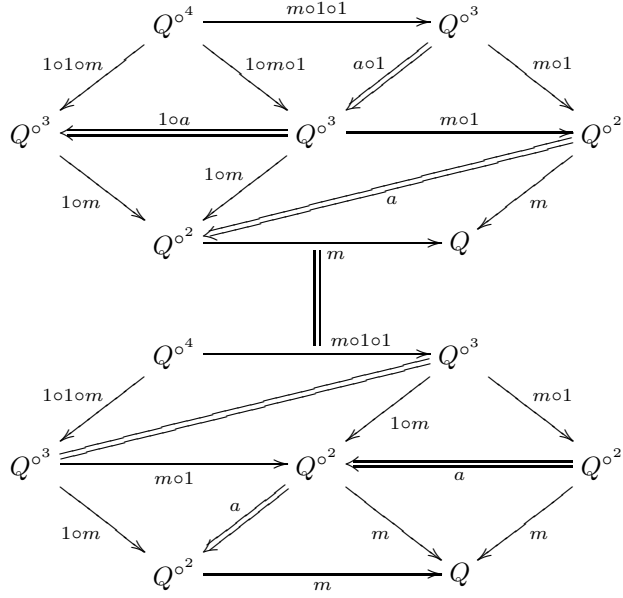

Finally, over $X^{[5]}$ we can define a bundle gerbe 
$(Q^{\circ^4},Y^{\circ^4},X^{[5]})$ 
where $Q^{\circ^4}$ and $Y^{\circ^4}$ 
are defined in the obvious way.  So for example, $Y^{\circ^4}$  
is the local-section-admitting surjection 
on $X^{[5]}$ with fiber 
$$
Y_{(x_4,x_5)}\times Y_{(x_3,x_4)}\times Y_{(x_2,x_3)}
\times Y_{(x_1,x_2)} 
$$
at a point $(x_1,x_2,x_3,x_4,x_5)\in X^{[5]}$.  Now the 
bundle gerbe morphism $m$ gives rise to five bundle 
gerbe morphisms $M_i\colon Q^{\circ^4}\to 
Q$, $i=1,\ldots 5$ covering the map $X^{[5]}\to X^{[2]}$ 
which sends $(x_1,x_2,x_3,x_4,x_5)$ to $(x_1,x_5)$.  The 
bundle gerbe morphisms $M_i$ are given as follows:  
$M_1 = m (m\circ 1)(m\circ 1\circ 1)$,  
$M_2 = m(m\circ 1)(1\circ m\circ 1)$, 
$M_3 = m(1\circ m)(1\circ m\circ 1)$,  
$M_4 = m(1\circ m)(1\circ 1\circ m)$,  
$M_5 = m(m\circ 1)(1\circ 1\circ m)$.   
Here we have abused notation and denoted for example 
by $m\circ 1$ the bundle gerbe morphism 
$Q^{\circ^3} \to Q^{\circ^2}$ which sends $u_{34}
\otimes u_{23}\otimes u_{12}$ to $\hat{m}(u_{34}\otimes 
u_{23})\otimes u_{12}$.    
Notice that $M_5$ can also be written as 
$M_5 = m(1\circ m)(m\circ 1\circ 1)$.  
It is not too hard to see that we have the following isomorphisms 
of $\cstar$ bundles on $X^{[5]}$.  We have $D_{M_1,M_2} 
= \pi_1^{-1}A$, $D_{M_2,M_3} = \pi_3^{-1}A$, $D_{M_3,M_4} 
= \pi_5^{-1}A$, $D_{M_4,M_5} = \pi_2^{-1}A^*$ and 
$D_{M_5,M_1} = \pi_4^{-1}A^*$.  From 
Lemma~\ref{lemma:bg transformations} 
there is an isomorphism 
$$
D_{M_1,M_2}\otimes D_{M_2,M_3}\otimes D_{M_3,M_4}
\otimes D_{M_4,M_5}\otimes D_{M_5,M_1} = D_{M_1,M_1},  
$$
and therefore , since $D_{M_1,M_1}$ is canonically trivialised, 
the $\cstar$ bundle $\d(A)$ on $X^{[5]}$ must be 
canonically trivialised.  Here $\d(A)$ is the $\cstar$ 
bundle given by 
$$
\d(A) = \pi_1^{-1}A\otimes \pi_2^{-1}A^*\otimes \pi_3^{-1}A
\otimes \pi_4^{-1}A^*\otimes \pi_5^{-1}A.  
$$
We finally require that the induced section $\d(a) = 
\pi_1^{-1}a\otimes \pi_2^{-1}a^*\otimes\pi_3^{-1}a \otimes 
\pi_4^{-1}a^*\otimes \pi_5^{-1}a$ of $\d(A)$ matches 
this canonical trivialisation.  This coherency condition 
on the section $a$ should actually be viewed as an equality 
of transformations of bundle gerbe morphisms as indicated in 
Figure~\ref{fig:associativity coherence for bundle gerbes}.    
Notice that this bit of theory is possible precisely 
because the $\pi_i$ are the face operators for a 
simplicial manifold.  All that we have said applies 
equally well to an arbitrary simplicial manifold.  
Hence we make the following definition.  

\begin{definition}[\cite{Ste}]  
\label{def:simplicial bg} 
A \emph{simplicial bundle gerbe} on a simplicial 
manifold $X_{\bullet} = \{X_p\}$ consists of the 
following data.  

\begin{enumerate} 

\item A bundle gerbe $(Q,Y,X_1)$ on $X_1$.  

\item A bundle gerbe morphism $m\colon d_0^{-1} 
Q\otimes d_2^{-1}Q \to d_1^{-1}Q$ over $X_2$.  

\item A transformation $a\colon m_1 \Rightarrow 
m_2$ between the two induced bundle gerbe morphisms 
$m_1$ and $m_2$ over $X_3$.  $m_1$ 
and $m_2$ are defined as in the following diagram.  
$$
\xymatrix{ 
d_0^{-1}(d_0^{-1}Q\otimes d_2^{-1}Q)\otimes 
d_2^{-1}d_2^{-1}Q \ar[d]_-{d_0^{-1}m
\otimes d_2^{-1}d_2^{-1}1_Q} \ar @2{-}[r] & 
d_1^{-1}d_0^{-1}Q\otimes d_3^{-1}(d_0^{-1}
Q\otimes d_2^{-1}Q) \ar[d]^-{d_1^{-1}d_0^{-1}
1_Q \otimes d_3^{-1}m}                  \\ 
d_0^{-1}d_1^{-1}Q\otimes d_2^{-1}d_2^{-1}Q 
\ar @2{-}[d] \ar @<-4.2ex> @2{->}[r]^-{a} & 
d_1^{-1}d_0^{-1}Q\otimes d_3^{-1}d_1^{-1}Q 
\ar @2{-}[d]                                   \\ 
d_2^{-1}(d_0^{-1}Q\otimes d_2^{-1}Q)  
\ar[d]_-{d_2^{-1}m} & d_1^{-1}(d_0^{-1}Q\otimes 
d_2^{-1}Q)   \ar[d]^-{d_1^{-1}m}             \\ 
d_2^{-1}d_1^{-1}Q \ar @2{-}[r] & d_1^{-1}d_1^{-1}Q.  } 
$$
So $m_1 = d_2^{-1}m\circ (d_0^{-1}m
\otimes d_2^{-1}d_2^{-1}1_Q)$ and $m_2 = 
d_1^{-1}m\circ (d_1^{-1}d_0^{-1}1_Q \otimes 
d_3^{-1}m)$.  Thus $a$ is a section of the 
$\cstar$ bundle $A = D_{m_1 ,m_2}$ 
over $X_3$.  

\item The transformation $a$ satisfies the coherency 
condition 
$$
d_0^{-1}a\otimes d_1^{-1}a^* \otimes d_2^{-1}a\otimes 
d_3^{-1}a^*\otimes d_4^{-1}a = 1, 
$$
where $1$ is the canonical section of the $\cstar$ 
bundle $\d(A)$ over $X_4$.  

\end{enumerate} 
\end{definition} 

Note that the coherency condition on the 
transformation $a$ can also be viewed as the 
commutativity of a diagram of the form 
Figure~\ref{fig:associativity coherence for bundle gerbes}.   
Clearly the notion of a simplicial bundle 
gerbe is a special case of Brylinski and 
McLaughlins definition of a simplicial 
gerbe \cite{BryMac}.  To recover the definition 
of simplicial gerbe from Definition~\ref{def:simplicial bg} 
above, simply replace each occurrence of the 
word `bundle gerbe' by the word `gerbe', 
`bundle gerbe morphism' by `gerbe morphism' 
and so on (strictly speaking we should insert 
certain canonical equivalences of gerbes where we 
have equalities of bundle gerbes, but this is of no real 
importance).  Note that the associator transformation 
of gerbe morphisms in the definition of a simplicial 
gerbe can be interpreted as a section of a certain 
line bundle on $X_3$, and the coherency condition 
on the transformation can be interpreted as a 
coherency condition on sections of line 
bundles on $X_4$, as above. 
  
We define a bundle $2$-gerbe to be a special case of 
the above definition.  

\begin{definition}[\cite{Ste}]  
\label{def:b2g} 
A \emph{bundle $2$-gerbe} consists of a quadruple of smooth 
manifolds $(Q,Y,X,M)$ where $\pi\colon X\to M$ 
is a smooth surjection admitting local sections 
and where $(Q,Y,X^{[2]})$ is a simplicial bundle 
gerbe on the simplicial manifold $X_{\bullet} 
= \{X_p\}$ with $X_p = X^{[p+1]}$ associated 
to $\pi\colon X\to M$.  
\end{definition} 

So given a bundle $2$-gerbe $(Q,Y,X,M)$, we have 
a bundle gerbe $(Q,Y,X^{[2]})$ and a bundle gerbe 
morphism $m\colon \pi_1^{-1}Q\otimes \pi_3^{-1}
Q\to \pi_2^{-1}Q$.  The bundle gerbe morphism 
$m$ consists of a pair of maps $(\hat{m},m)$, 
where $m\colon Y_{123} = \pi_1^{-1}Y\times_{X^{[3]}}\pi_3^{-1}Y 
\to Y_{13} = \pi_2^{-1}Y$ is a map commuting with the projections to 
$X^{[3]}$ and $\hat{m}\colon \pi_1^{-1}Q\otimes 
\pi_3^{-1}Q\to \pi_2^{-1}Q$ covers $m^{[2]}\colon 
Y_{123}^{[2]}\to 
Y_{13}^{[2]}$ and commutes with the bundle gerbe 
products on $\pi_1^{-1}Q\otimes \pi_3^{-1}Q$ and 
$\pi_2^{-1}Q$.  So fiberwise $m$ is a map 
$$
m\colon Y_{(x_2 ,x_3)} \times Y_{(x_1 ,x_2)} \to 
Y_{(x_1 ,x_3)} 
$$
for $(x_1 ,x_2 ,x_3) \in X^{[3]}$ and $\hat{m}$ 
is a map 
$$
\hat{m} \colon Q_{(y_{23},y_{23}^{'})}\otimes 
Q_{(y_{12},y_{12}^{'})} \to Q_{(m(y_{23},y_{12}),
m(y_{23}^{'},y_{12}^{'}))} 
$$
for $(y_{23},y_{12}),\ (y_{23}^{'},y_{12}^{'}) 
\in Y_{123}$.  Thus for each pair of points  
$(x_1 ,x_2)$ lying in the same fiber of $\pi 
\colon X\to M$, we obtain a $\cstar$ groupoid 
$\textbf{Gr}(Q)_{(x_1 ,x_2)}$.  Given a triple 
of points $(x_1 ,x_2 ,x_3)$ lying in the 
same fiber of $\pi\colon X\to M$ the bundle 
gerbe morphism $m$ gives rise to a 
functor $m\colon \textbf{Gr}(Q)_{(x_2 ,x_3)} \times 
\textbf{Gr}(Q)_{(x_1 ,x_2)}\to \textbf{Gr}(Q)_{(x_1 ,
x_3)}$ as explained in Section~\ref{sec:four}.  
Let us denote the action of the functor $m$ on a 
pair of objects $(y_{23},y_{12})$ of $\textbf{Gr}
(Q)_{(x_2,x_3)}\times \textbf{Gr}(Q)_{(x_1,x_2)}$ 
by $y_{23}\circ y_{12}$.     
The transformation $a$ gives rise to a natural 
transformation, also denoted $a$, between the 
functors bounding the following diagram.  
$$
\xymatrix{ 
\textbf{Gr}(Q)_{(x_3,x_4)}\times \textbf{Gr}
(Q)_{(x_2,x_3)}\times \textbf{Gr}(Q)_{(x_1,x_2)} 
\ar[r]^-{m\times 1} \ar[d]_-{1\times m} & \textbf{Gr}(Q)_{(x_2,
x_4)}\times \textbf{Gr}(Q)_{(x_1,x_2)} \ar[d]^-{m} 
\ar @2{->}[dl]^-{a}                                 \\ 
\textbf{Gr}(Q)_{(x_3,x_4)}\times \textbf{Gr}(Q)_{
(x_1,x_3)} \ar[r]_-{m} & \textbf{Gr}(Q)_{(x_1,x_4)} } 
$$
The coherency condition on the transformation $a$ 
of bundle gerbe morphisms can be  
viewed as an \emph{associativity coherence} condition 
on the natural transformation $a$.  
Let us briefly recall the definition of a 
\emph{bicategory} \cite{Ben}.  A bicategory 
$\mathcal{B}$ consists of objects $A,B,C,\ldots$ 
and for each pair of objects $A$ and $B$ a category 
$\Hom(A,B)$.  The objects of $\Hom(A,B)$ 
are called \emph{$1$-arrows} or \emph{$1$-cells} of $\mathcal{B}$ 
and the arrows of $\Hom(A,B)$ are called \emph{$2$-arrows} 
or \emph{$2$-cells} of $\mathcal{B}$.  A $2$-cell $\phi$ 
between $1$-cells $\a$ and $\b$ of $\Hom
(A,B)$ is denoted $\phi\colon \a\Rightarrow \b$.  
Given three objects 
$A$, $B$ and $C$ of $\mathcal{B}$ there is a 
composition functor 
$\Hom(B,C)\times \Hom(A,B) \to \Hom
(A,C)$ whose action on a pair of objects $(\a,\b)$ 
of $\Hom(B,C)\times \Hom(A,B)$ 
($1$-cells of $\mathcal{B}$) is 
denoted $\a\circ \b$ and similarly for $2$-cells.  
The composition functor is associative up to 
a coherent isomorphism.  This means that given 
objects $A$, $B$, $C$ and $D$ of $\mathcal{B}$ 
with $1$-cells $\a \in \Hom(A,B)$, 
$\b \in \Hom(B,C)$ and $\c \in \Hom
(C,D)$ then there is an isomorphism 
$$
a(\c,\b,\a)\colon (\c\circ \b)\circ \a 
\Rightarrow \c \circ (\b\circ \a) 
$$
in $\Hom(A,C)$ which is natural in $A$, $B$ 
and $C$.  The natural isomorphism $a$ is called 
the \emph{associator} natural isomorphism.  
The associativity coherence condition 
means that the well known pentagonal diagram  
commutes.  One also requires that 
for every object $A$ of $\mathcal{B}$ there is a 
$1$-arrow $1_A$ of $\Hom(A,A)$ and for every 
$1$-arrow $\a \in \Hom(A,B)$ of $\mathcal{B}$ 
there are left and right identity isomorphisms 
$L_{\a}\colon \a\circ 1_A \Rightarrow \a$ and 
$R_{\a}\colon 1_B \circ \a \Rightarrow \a$ which are 
natural in the $1$-arrows $\a$.  These isomorphisms 
are finally required to satisfy the coherency 
condition that the following diagram commutes.  
$$
\xymatrix{ 
(\b\circ 1_B)\circ \a \ar @2{->}[d]_-{L_{\b}\circ 
1_{\a}} \ar @2{->}[r]^-{a(\a,1_B,\b)} & 
\b\circ (1_B\circ \a) \ar @2{->}[d]^-{1_{\b}\circ R_{\a}} \\ 
\b \circ \a \ar @2{-}[r] & \b\circ \a                     } 
$$
where $\a$ is a $1$-arrow of $\Hom(A,B)$ and 
$\b$ is a $1$-arrow of $\Hom(B,C)$.  A bicategory 
in which all of the natural isomorphisms $a$, $L$ 
and $R$ are the identities is a $2$-category.  

One can define the notion of a \emph{biequivalence} 
between bicategories; we will refer to \cite{Ben} 
for this.  One can show \cite{GorPowStr} that every 
bicategory is biequivalent to a $2$-category.  We also 
have the notion of a bigroupoid.  

\begin{definition}[\cite{Bre}] 
\label{def:bigroupoid} 
A \emph{bigroupoid} consists of a bicategory $\mathcal{B}$ 
which satisfies the following two additional axioms.  

\begin{enumerate} 
\item $1$-arrows are coherently invertible.  This means that 
if $\a$ is a $1$-arrow of $\Hom(A,B)$ then there is a 
$1$-arrow $\b$ of $\Hom(B,A)$ together with 
$2$-arrows $\phi\colon \b\circ \a\Rightarrow 1_A$ 
of $\Hom(A,A)$ and 
$\psi\colon \a\circ \b\Rightarrow 1_B$ in $\Hom(B,B)$.  

\item All $2$-arrows are invertible. 
\end{enumerate} 
By a $\cstar$ bigroupoid we mean a bigroupoid $\mathcal{B}$ 
in which the automorphism group of every $1$-arrow is isomorphic 
to $\cstar$.   
\end{definition} 
  
We have the following Proposition.  
\begin{proposition}[\cite{Ste}] 
For each point $m$ of $M$, the restriction of a bundle 
$2$-gerbe $(Q,Y,X,M)$ to the point $m$ gives rise to a 
family of $\cstar$ bigroupoids $\mathcal{Q}_m$. 
\end{proposition} 

We take as the objects of $\mathcal{Q}_m$ the points of 
$X_m = \pi^{-1}(m)$.  Given two such points $x_1$ and 
$x_2$ we define the category $\Hom(x_1,x_2)$ 
to be the category $\textbf{Gr}(Q)_{(x_1,x_2)}$ defined 
above.  It is clear that the bundle gerbe morphism 
$m$ provides the composition functor and that 
the transformation $a$ plays the role of the associator 
natural isomorphism.  All we have to do then is to define 
left and right identity morphisms and show that they 
are compatible with $a$.  We will not do this here and refer 
instead to \cite{Ste}.  Thus we can think of a bundle $2$-gerbe 
as being a `bundle of bigroupoids'.     

\section{The Homotopy Bigroupoid and the Tautological 
Bundle 2-Gerbe.} 
\label{sec:six} 

An important example of a bigroupoid is the so-called 
\emph{homotopy bigroupoid} or \emph{fundamental 
bigroupoid} $\mathbf{\Pi}_2 (X)$ 
associated to a topological 
space $X$ (see \cite{Bat}).  $\mathbf{\Pi}_2(X)$ is defined as follows.  
The objects of $\mathbf{\Pi}_2(X)$ are the points $x$ of 
$X$.  Given two points $x_1$ and $x_2$ the 
category $\Hom(x_1,x_2)$ is defined 
to have as objects ($1$-cells of $\mathbf{\Pi}_2(X)$) 
the paths $\c\colon I\to X$ with $\c(0) = x_1$ 
and $\c(1) = x_2$ where $I$ denotes the unit interval 
$[0,1]$.  Given two such paths $\c_1$ and $\c_2$ 
the set of $2$-cells $\c_1 \Rightarrow \c_2$  is 
defined to be the set of homotopy classes $[\mu]$  
of maps $\mu\colon I\times I\to X$ such that $\mu(0,t) 
= \c_1(t)$, $\mu(1,t) = \c_2(t)$, $\mu(s,0) = 
x_1$ and $\mu(s,1) = x_2$.  Two such maps $\mu$ and 
$\mu'$ belong to the same homotopy class if there 
is a map $H\colon I\times I\times I\to X$ such that 
$H(0,s,t) = \mu(s,t)$, $H(1,s,t) = \mu'(s,t)$,  
$H(r,0,t) = \c_1(t)$, $H(r,1,t) = \c_2(t)$, 
$H(r,s,0) = x_1$ and $H(r,s,1) = x_2$.  To define the 
composite $2$-cell $[\lambda] [\mu]\colon \c_1\Rightarrow 
\c_3$ for $2$-cells $[\mu]\colon \c_1\Rightarrow 
\c_2$ and $[\lambda]\colon \c_2\Rightarrow \c_3$ 
we choose representatives $\mu$ and $\lambda$ of 
$[\mu]$ and $[\lambda]$ respectively and define 
$[\lambda][\mu]$ to be the homotopy class of the 
map 
$$
(\lambda \mu)(s,t) = \begin{cases} 
                     \mu (2s,t)\ s\in [0,\frac{1}{2}],\ 
                     t \in [0,1],                         \\ 
                     \lambda(2s-1,t),\ s\in [\frac{1}{2},1],\ 
                     t\in [0,1].                          \\ 
                      \end{cases} 
$$
It is straightforward to check that this law 
of composition is well defined and is associative.  Notice that 
every $2$-cell of $\mathbf{\Pi}_2(X)$ is invertible.  We need to 
define the composition functor 
$$
m\colon\Hom(x_2,x_3)\times \Hom(x_1,x_2)\to 
\Hom(x_1,x_3).  
$$
If $\c_{23}$ is a $1$-arrow of $\Hom(x_2,x_3)$ 
and $\c_{12}$ is a $1$-arrow of $\Hom(x_1,x_2)$ 
then we define $m(\c_{23},\c_{12})$ to be the 
path $\c_{23}\circ \c_{12}\colon I\to X$ given by 
$$
(\c_{23}\circ \c_{12})(t) = \begin{cases} 
                            \c_{12}(2t),\ t\in [0,\frac{1}{2}], \\ 
                            \c_{23}(2t-1),\ t\in [\frac{1}{2},1]. \\ 
                            \end{cases} 
$$
If $[\mu_{23}]\colon \c_{23}\Rightarrow 
\c'_{23}$ in $\Hom(x_2,x_3)$  
and $[\mu_{12}]\colon \c_{12}\Rightarrow \c'_{12}$  
in $\Hom
(x_1,x_2)$ are $2$-arrows, we define $m([\mu_{23}],
[\mu_{12}]) = [\mu_{23}\circ \mu_{12}]$ to 
be the homotopy class of the map 
$$
(\mu_{23}\circ \mu_{12})(s,t) = \begin{cases} 
                              \mu_{12}(s,2t),\ s\in [0,1],\ t\in 
                              [0,\frac{1}{2}],                    \\ 
                               \mu_{23}(s,2t-1),\ s\in [0,1],\ t\in 
                              [\frac{1}{2},1],                      \\ 
                              \end{cases} 
$$
where $\mu_{23}$ is a representative of the 
homotopy class $[\mu_{23}]$ and $\mu_{12}$ is 
a representative of the homotopy class $[\mu_{12}]$.  
Note that the map $\mu_{23}\circ \mu_{12}\colon 
I\times I\to X$ defines a homotopy with endpoints fixed 
between $\c_{23}\circ \c_{12}$ and $\c'_{23}\circ \c'_{12}$.    
It is straightforward to check that this defines 
a functor.  We now need to define identity $1$-arrows 
and identity $2$-arrows.  Given an object $x$ of $\mathbf{\Pi}_2 (X)$, 
we define $1_x$ to be the constant path at $x$ and 
the identity $2$-arrow $1_x \Rightarrow 1_x$ to be the 
constant homotopy from the constant path to itself.  

Next we define the associator isomorphism.  Given $1$-arrows 
$\c_{34}$ in $\Hom(x_3,x_4)$, $\c_{23}$ in 
$\Hom(x_2,x_3)$ and $\c_{12}$ in $\Hom(x_1,x_2)$ 
we need to define a $2$-arrow 
$$
a(\c_{34},\c_{23},\c_{12})\colon 
(\c_{34}\circ \c_{23})\circ \c_{12}\Rightarrow \c_{34}\circ 
(\c_{23}\circ \c_{12}).
$$
There is a standard choice for 
$a(\c_{34},\c_{23},\c_{12})$ --- see for example \cite{Spa}.  
We set $a(\c_{34},\c_{23},\c_{12})$ equal to the homotopy 
class of the map $\bar{a}(\c_{34},\c_{23},\c_{12})\colon 
I\times I\to X$ given by 
\begin{equation} 
\label{eq:associator section} 
\bar{a}(\c_{34},\c_{23},\c_{12})(s,t) = \begin{cases} 
                      \c_{12}(\frac{4t}{2-s}),\ s\in [0,1],\ 
                      t\in [0,\frac{2-s}{4}],                  \\ 
                      \c_{23}(4t-2+s),\ s\in [0,1],\ 
                      t\in [\frac{2-s}{4},\frac{3-s}{4}],       \\ 
                      \c_{34}(\frac{4t-3+s}{1+s}),\ s\in [0,1],\ 
                      t\in [\frac{3-s}{4},1].                  
                      \end{cases} 
\end{equation} 
One can check that the assignment of the $2$-arrow 
$a(\c_{34},\c_{23},\c_{12})$ of $\Hom(x_1,x_4)$ 
to the $1$-arrow $(\c_{34},\c_{23},\c_{12})$ of 
$\Hom(x_3,x_4)\times \Hom(x_2,x_3)\times 
\Hom(x_1,x_2)$ is a natural transformation 
$m\circ (m\times 1)\Rightarrow m\circ (1\times m)$.  
We now have to check that the natural transformation 
$a$ satisfies the associativity coherence condition.  
This means that we have to check that the diagram 
of $2$-arrows in Figure~\ref{fig:associativity 
coherence for bicategories} is the identity $2$-arrow from $((\c_{45}\circ 
\c_{34})\circ \c_{23})\circ \c_{12}$ to itself.  
\begin{figure} 
$$ 
\xymatrix@!C=40pt{ 
& & (\c_{45}\circ (\c_{34}\circ \c_{23}))\circ \c_{12} 
\ar @2{->}[drr]^-{a} & & \\ 
((\c_{45}\circ \c_{34})\circ \c_{23})\circ \c_{12} 
\ar @2{->}[urr]^-{a\circ 1_{\c_{12}}\ \ \ \ \ } 
& & & & \c_{45}\circ ((\c_{34}\circ \c_{23})\circ \c_{12}) 
\ar @2{->}[dl]_-{1_{\c_{45}}\circ a}   \\ 
& (\c_{45}\circ \c_{34})\circ (\c_{23}\circ \c_{12}) 
\ar @2{->}[ul]_-{a^{-1}} 
& & \c_{45}\circ (\c_{34}\circ (\c_{23}\circ \c_{12})) 
\ar @2{->}[ll]^-{a^{-1}}     } 
$$ 
\caption{Associativity coherence condition for $\Pi_2(X)$} 
\label{fig:associativity coherence for bicategories} 
\end{figure} 
We will omit the proof of this fact and refer to 
\cite{Ste} where an explicit homotopy between the 
composed $2$-arrow from $((\c_{45}\circ \c_{34})\circ 
\c_{23})\circ \c_{12}$ to itself and the identity $2$-arrow 
is given.  To show that $\mathbf{\Pi}_2(X)$ is a bicategory, 
we need to produce left and right identity isomorphisms.  
If $\c\in \Hom(x_1,x_2)$ then $L(\c)$ is a 
$2$-arrow $\c\Rightarrow \c\circ 1_{x_1}$.  We define 
$L(\c)$ to be the homotopy class of the map 
$$
(s,t)\mapsto \begin{cases} 
             x_1,\ t\in [0,\frac{s}{2}], \\ 
             \c(\frac{2t-s}{2-s}),\ t\in [\frac{s}{2},1]. 
             \end{cases} 
$$
Similarly if $\c\in \Hom(x_1,x_2)$ then $R(\c)$ 
is a $2$-arrow $R(\c)\colon \c \Rightarrow 1_{x_2}\circ \c$.  
We set $R(\c)$ equal to the homotopy class of the map 
$$
(s,t)\mapsto \begin{cases} 
             \c((s+1)t),\ t\in [0,\frac{1}{s+1}], \\ 
             x_2,\ t\in [\frac{1}{s+1},1].         
             \end{cases} 
$$
One can check (see \cite{Ste}) that the assignments 
$\c \mapsto L(\c)$ and $\c\mapsto R(\c)$ define natural 
transformations and that moreover these natural transformations 
are compatible with $a$.  Hence $\mathbf{\Pi}_2(X)$ is an example of 
a bicategory.  One can also show that the    
$1$-arrows of $\mathbf{\Pi}_2(X)$ are 
coherently invertible and,       
as mentioned 
earlier, all $2$-arrows of $\mathbf{\Pi}_2(X)$ are invertible.  
Therefore $\mathbf{\Pi}_2(X)$ is a 
bigroupoid - the homotopy bigroupoid 
of $X$.  

We will now use this description of the 
homotopy bigroupoid $\mathbf{\Pi}_2(X)$ of $X$ to 
define the tautological bundle $2$-gerbe of \cite{CarMurWan}  
over a $3$-connected manifold $M$.  Recall that we start 
with a closed four form $\Theta$ on $M$ 
with integral periods, representing a class in 
$H^4 (M;\Z)$.  We then form the path fibration 
$\pi\colon \mathcal{P}M\to M$, where $\mathcal{P}M$ 
is the Frechet manifold consisting of 
piecewise smooth paths $\c\colon [0,1]\to M$, 
$\c(0) = m_0$ where $m_0$ is a basepoint of $M$, and 
where $\pi$ is the map sending such a path $\c$ to its 
endpoint $\c(1)$.  The fibration $\pi\colon \mathcal{P}
M\to M$ has fiber $F$ equal to the space of piecewise smooth 
loops in $M$, $\Omega M$.    

We will define a simplicial bundle gerbe 
on the simplicial manifold $X_{\bullet} = \{X_p\}$ 
with $X_p = X^p$ and with 
face and degeneracy operators $d_i\colon 
X^{p+1}\to X^p$, $s_i\colon 
X^p \to X^{p+1}$ given respectively 
by 
\begin{eqnarray*} 
d_i (x_1,\ldots,x_{p+1}) & = & (x_1,\ldots,x_{i-1},x_{i+1},
\ldots,x_{p+1})                                             \\ 
s_i (x_1,\ldots,x_p) & = & (x_1,\ldots,x_i,x_i,\ldots,
x_p). 
\end{eqnarray*} 
Performing the construction for the simplicial 
manifold $X_{\bullet}$ above with $X = \Omega M$ fiber by 
fiber on $\mathcal{P}M$ will define the tautological 
bundle $2$-gerbe.  We start with a $2$-connected manifold $X$ and a 
closed $3$-form $\omega$ on $X$ with integral periods.  
We construct a bundle gerbe on $X^2 = X\times X$ in 
the usual way.  We define a fibering $Y\to X^2$ with 
fiber $Y_{(x_1,x_2)}$ at $(x_1,x_2)\in X^2$ equal to 
the space of piecewise smooth paths $\a\colon I\to 
X$ with $\a(0) = x_1$ and $\a(1) = x_2$.  Next we define 
a $\cstar$ bundle $Q\to Y^{[2]}$ whose fiber at $(\a,\b)
\in Y^{[2]}$ is all equivalence classes $[\mu,z]$ where 
$z\in \cstar$ and $\mu\colon I^2  
\to X$ is a homotopy with 
endpoints fixed between $\a$ and $\b$, that is $\mu(0,t) = 
\a(t)$, $\mu(1,t) = \b(t)$, $\mu(s,0) = x_1$ and 
$\mu(s,1) = x_2$.  The equivalence 
relation $\equivalencerelation$ is defined by declaring   
$(\mu_1,z_1)\equivalencerelation (\mu_2,z_2)$ if 
for any homotopy $F\colon I^3 \to 
X$ with endpoints fixed between $\mu_1$ and $\mu_2$ 
we have 
$$
z_2 = z_1 \exp (\int_{I^3} F^* (\omega)).  
$$
Here we say that $F$ is a homotopy with endpoints fixed between 
$\mu_1$ and $\mu_2$ if we have $F(0,s,t) = \mu_1(s,t)$, 
$F(1,s,t) = \mu_2(s,t)$, $F(r,0,t)= \a(t)$, 
$F(r,1,t) = \b(t)$, $F(r,s,0)= x_1$ and 
$F(r,s,1)= x_1$.  One can define an associative 
product $m_Q$ on $Q\to Y^{[2]}$ as in \cite{CarMurWan} 
by setting 
$m_Q ([\mu,z]\otimes [\nu,w]) = [\mu \nu,zw]$, where 
$\mu \nu\colon I^2\to X$ is defined by 
$$
(\mu \nu)(s,t) = \begin{cases} 
                      \nu(2s,t),\ s\in [0,\frac{1}{2}],\ 
                      t\in [0,1],                       \\ 
                      \mu(2s-1,t),\ s\in [\frac{1}{2},1],\ 
                      t\in [0,1].   
                      \end{cases} 
$$
One can check, see \cite{CarMurWan}, that this 
is well defined and associative.  Next we define 
a bundle gerbe morphism $m\colon d_0^{-1}Q
\otimes d_2^{-1}Q\to d_1^{-1}Q$ with $m = 
(\hat{m},m)$.  So fiberwise $m$ will be a map 
$Y_{(x_2,x_3)}\times Y_{(x_1,x_2)}\to Y_{(x_1,x_3)}$.  
$m$ is defined by the composition functor in 
the bigroupoid $\mathbf{\Pi}_2(X)$ so $m(\a,\b) = \a\circ \b$ where 
$\a\circ \b\colon I\to X$ is the path from $x_1$ 
to $x_3$ given by 
$$
(\a\circ \b)(t) = \begin{cases} 
                  \b(2t),\ t\in [0,\frac{1}{2}], \\ 
                  \a(2t-1),\ t\in [\frac{1}{2},1]. 
                  \end{cases} 
$$
The map $\hat{m}\colon d_0^{-1}Q\otimes d_2^{-1}Q
\to d_1^{-1}Q$ covering $m^{[2]}$ is defined by 
$\hat{m}([\mu_{23},z_{23}]\otimes [\mu_{12},z_{12}]) 
= [\mu_{23}\circ \mu_{12},z_{23}z_{12}]$, where 
$\mu_{23}\circ \mu_{12}$ is defined by the action of 
the composition functor $m$ in the bigroupoid $\mathbf{\Pi}_2(X)$ 
on $2$-arrows.  Hence $\mu_{23}\circ \mu_{12}\colon 
I^2\to X$ is the homotopy given by 
$$
(\mu_{23}\circ \mu_{12})(s,t) = \begin{cases} 
                                \mu_{12}(s,2t),\ s\in [0,1],\ 
                                t\in [0,\frac{1}{2}],         \\ 
                                \mu_{23}(s,2t-1),\ s\in [0,1],\ 
                                t\in [\frac{1}{2},1].  
                                \end{cases} 
$$
Again, one can check (see \cite{CarMurWan}), that this 
is well defined and commutes with the bundle gerbe 
products.  As usual, $m$ defines two bundle 
gerbe morphisms $m_1 = (\hat{m}_1,m_1)$, 
$m_2 = (\hat{m}_2,m_2)$ between the appropriately 
defined bundle gerbes on $X^4$.  So fiberwise $m_1$ and 
$m_2$ are maps $Y_{(x_3,x_4)}\times Y_{(x_2,x_3)}\times 
Y_{(x_1,x_2)}\to Y_{(x_1,x_4)}$ which are given by 
$m_1(\a_{34},\a_{23},\a_{12}) = (\a_{34}\circ \a_{23})
\circ \a_{12}$ , $m_2(\a_{34},\a_{23},\a_{12}) = 
\a_{34}\circ (\a_{23}\circ \a_{12})$.  $\hat{m}_1$ 
and $\hat{m}_2$ are defined in an analogous fashion.  
As we have already seen, there is a homotopy $m_1 \simeq 
m_2$.  We can use this homotopy to write down a section 
$\hat{a}$ which trivialises the $\cstar$ bundle $(m_1,
m_2)^{-1}Q$ on $Y\circ Y\circ Y$.  We have $\hat{a}(\a_{34},
\a_{23},\a_{12}) = [\bar{a}(\a_{34},\a_{23},\a_{12}),1]$ 
where $\bar{a}(\a_{34},\a_{23},\a_{12})\colon I^2 \to X$ 
is defined in equation~\ref{eq:associator section}.  
Recall that the associator 
natural isomorphism for the bigroupoid $\mathbf{\Pi}_2(X)$ 
is defined via $\bar{a}$.  The fact that this is a natural 
isomorphism is exactly the requirement that $\hat{a}$ 
descends to a section $a$ of the bundle $A$ on $X^4$.  
Finally, one needs to show that $a$ satisfies the 
coherency condition over $X^5$ or, alternatively, that 
$\hat{a}$ satisfies the analogous coherency condition.  
Let $\d(\hat{a})$ denote the $2$-arrow in 
Figure~\ref{fig:associativity coherence for bicategories} from 
$(\c_{45}\circ (\c_{34}\circ \c_{23}))\circ \c_{12}$ 
to itself.  In \cite{Ste} an explicit homotopy from 
$\d(\hat{a})$ to the identity $2$-arrow at $(\c_{45}\circ 
(\c_{34}\circ \c_{23}))\circ \c_{12}$ was written down.  
One checks easily that the pullback of $\omega$ by 
this homotopy is zero.  This shows that $a$ satisfies the 
required coherency condition.        

\section{The \v{C}ech $3$-class associated to a Bundle $2$-gerbe.} 
\label{sec:seven} 

Let $(Q,Y,X,M)$ be a bundle $2$-gerbe.  We will 
explain how to construct a $\cstar$ valued \v{C}ech 
$3$-cocycle associated to $Q$.  Choose an open covering 
$\{U_i\}_{i\in I}$ of $M$ all of whose finite intersections 
are empty or contractible and such that there exist 
local sections $s_i\colon U_i\to X$ of $\pi\colon X\to M$.  
Form maps $(s_i,s_j)\colon U_{ij}\to X^{[2]}$ by sending 
a point $m$ of $U_{ij}$ to the point $(s_i(m),s_j(m))$ 
of $X^{[2]}$.  Let $(Q_{ij},Y_{ij},U_{ij})$ denote 
the pullback of the bundle gerbe $(Q,Y,X^{[2]})$ 
to $U_{ij}$ via $(s_i,s_j)$.  Therefore $Y_{ij}\to U_{ij}$ 
is a local-section-admitting surjection and the 
fiber $(Y_{ij})_m$ of $Y_{ij}$ at $m\in U_{ij}$ is 
$Y_{(s_i(m),s_j(m))}$.  

Since $U_{ij}$ is contractible, the bundle gerbe 
$(Q_{ij},Y_{ij},U_{ij})$ is trivial.  Hence there 
is a $\cstar$ bundle $P_{ij}$ on $Y_{ij}$ and an 
isomorphism $Q_{ij}\to \d(P_{ij})$ over $Y_{ij}^{[2]}$ 
which commutes with the bundle gerbe products on $Q_{ij}$ 
and the trivial bundle gerbe $\d(P_{ij})$.  The bundle 
gerbe morphism $m\colon \pi_1^{-1}Q\otimes 
\pi_3^{-1}Q\to \pi_2^{-1}Q$ pulls back to define 
a bundle gerbe morphism $Q_{jk}\otimes Q_{ij}\to Q_{ik}$, 
also denoted $m$.  In particular there is a map 
$m\colon Y_{jk}\times_{M}Y_{ij}\to Y_{ik}$ covering the 
identity on $U_{ijk}$.  Let $\hat{P}_{ijk} = 
P_{jk}\otimes m^{-1}P_{ik}^{*}\otimes P_{ij}$.  
Thus $\hat{P}_{ijk}$ is a $\cstar$ bundle on 
$Y_{jk}\times_{M}Y_{ij}$.  Note that there is an 
isomorphism $\d(\hat{P}_{ijk})\to Q_{jk}\otimes 
(m^{[2]})^{-1}Q_{ik}^{*}\otimes Q_{ij}$ which 
commutes with the respective bundle gerbe products.  
Moreover, the $\cstar$ bundle $Q_{jk}\otimes 
(m^{[2]})^{-1}Q_{ik}^*\otimes 
Q_{ij}$ has a canonical trivialisation provided by 
the bundle gerbe morphism $m$.   
The following Lemma follows easily from 
Example~\ref{ex:difference of two trivialisations}. 

\begin{lemma} 
\label{lemma:Cech 3-class} 
Suppose $(P,X,M)$ and $(Q,Y,M)$ are bundle 
gerbes with a bundle gerbe morphism $f
\colon P\to Q$.  If $P$ and $Q$ are both trivial, 
so there exist $\cstar$ bundles $T_P$ and $T_Q$ 
on $X$ and $Y$ respectively, with $\d(T_P) = P$ 
and $\d(T_Q) = Q$, then the bundle $T_P\otimes 
f^{-1}T_Q^*$ descends to $M$.              
\end{lemma} 
 
Applying this result we see that $\hat{P}_{ijk}$ descends 
to a $\cstar$ bundle $P_{ijk}$ on $U_{ijk}$.  Next, over 
$U_{ijkl}$ we have two induced bundle gerbe morphisms 
$m_1, m_2 \colon Q_{kl}\otimes Q_{jk}\otimes Q_{ij}
\to Q_{il}$.  By Lemma~\ref{lemma:bg transformations} 
the $\cstar$ bundle $(m_1,m_2)^{-1}Q_{il}$ on 
$Y_{ijkl} = Y_{kl}\times_M Y_{jk}\times_M Y_{ij}$ 
descends to  a $\cstar$ bundle $A_{ijkl}$ on $U_{ijkl}$, 
and it is clear that $A_{ijkl} = (s_i,
s_j,s_k,s_l)^{-1}A$.  We will show that there is an 
isomorphism 
$$
A_{ijkl} = P_{jkl}\otimes P_{ikl}^*\otimes P_{ijl}
\otimes P_{ijk}^* 
$$
of $\cstar$ bundles on $U_{ijkl}$.  Recall that the 
map $m_1\colon Y_{ijkl}\to Y_{il}$ is defined by composition: 
$Y_{ijkl}\stackrel{m\times 1}{\to}Y_{ijl}\stackrel{m}{\to} 
Y_{il}$, where $Y_{ijl} = Y_{jl}\times_M Y_{ij}$.  
It is not hard to show that  
$P_{kl}\otimes P_{jk}\otimes P_{ij}\otimes 
m_1^{-1}P_{il}^*   
 \simeq  \pi_{Y_{ijkl}}^{-1}(P_{ijl}\otimes P_{jkl})$.   
Similarly we get another isomorphism 
$P_{kl}\otimes P_{jk}\otimes P_{ij}\otimes m_2^{-1}P_{il}^* 
\simeq \pi_{Y_{ijkl}}^{-1}(P_{ijk}\otimes P_{ikl})$. 
Since we have an isomorphism $m_1^{-1}P_{il}^*\otimes 
m_2^{-1}P_{il} \simeq \pi_{Y_{ijkl}}^{-1}A_{ijkl}$, we 
get the required isomorphism 
$A_{ijkl} \simeq P_{jkl}\otimes P_{ikl}^*\otimes 
P_{ijl}\otimes P_{ijk}^*$  
over $U_{ijkl}$.  Now choose sections $\sigma_{ijk}$ 
of $P_{ijk}$ over $U_{ijk}$ and define $g_{ijkl}\colon 
U_{ijkl}\to \cstar$ by 
$$
\sigma_{jkl}\otimes \sigma_{ikl}^*\otimes \sigma_{ijl}
\otimes \sigma_{ijk}^*\cdot g_{ijkl} = a_{ijkl}. 
$$
One can show that $g_{ijkl}$ is a \v{C}ech 
$3$-cocycle.  We have the following Proposition. 
\begin{proposition}[\cite{Ste}]  
\label{prop:cech 3-class} 
$g_{ijkl}$ satisfies the \v{C}ech $3$-cocycle condition 
$$
g_{jklm}g_{iklm}^{-1}g_{ijlm}g_{ijkm}^{-1}g_{ijkl} = 1, 
$$
and hence is a representative of a class in 
$\check{H}^3 (M;\underline{\C}^{\times}_M) = H^4 (M;\Z)$. 
\end{proposition} 

There is another method of calculating the \v{C}ech 
$3$-cocycle $g_{ijkl}$ which is similar in spirit to 
the method used to calculate the \v{C}ech representative 
of the Dixmier-Douady class of a bundle gerbe.  Let 
$(Q,Y,X,M)$ be a bundle $2$-gerbe.  Choose an open cover 
$\{U_i\}_{i\in I}$ of $M$ all of whose finite non-empty 
intersections are contractible and such that there exist 
local sections $s_i\colon U_i\to X$ of $\pi$.  Form the 
maps $(s_i,s_j)\colon U_{ij}\to X^{[2]}$ as above 
and again denote the pullback of the bundle gerbe 
$(Q,Y,X^{[2]})$ to $U_{ij}$ via $(s_i,s_j)$ by 
$(Q_{ij},Y_{ij},U_{ij})$.  In certain circumstances, for 
instance if $\pi_Y\colon Y\to X^{[2]}$ is a 
fibration, one can choose sections $\sigma_{ij}\colon 
U_{ij}\to Y_{ij}$ of $\pi_{Y_{ij}}\colon Y_{ij}\to U_{ij}$.  
Note that in general one would only be able to 
choose an open cover $\{U_{ij}^{\a}\}_{\a\in\Sigma_{ij}}$ 
of $U_{ij}$ such that there were local sections 
$\sigma_{ij}^{\a}\colon U_{ij}^{\a}\to Y_{ij}$ of 
$\pi_{Y_{ij}}$.  We will assume here that we are
in the former situation described above.  
For ease of notation denote 
$m(\sigma_{jk},\sigma_{ij})$ by $\sigma_{jk}\circ 
\sigma_{ij}$.  Then we have a map $(\sigma_{ik},
\sigma_{jk}\circ \sigma_{ij})\colon U_{ijk}\to Y_{ik}
^{[2]}$ which sends $m\in U_{ijk}$ to $(\sigma_{ik}
(m),(\sigma_{jk}\circ \sigma_{ij})(m))\in Y_{ik}^{[2]}$.  
Let $Q_{ijk}$ denote the pullback $\cstar$ bundle 
$(\sigma_{ik},\sigma_{jk}\circ \sigma_{ij})^{-1}Q_{ik}$ 
on $U_{ijk}$.  We then have 
\begin{eqnarray*} 
(\sigma_{kl}\circ (\sigma_{jk}\circ \sigma_{ij}), 
\sigma_{il})^{-1}Q & = & (\sigma_{kl}\circ (\sigma_{jk}
\circ \sigma_{ij}),\sigma_{kl}\circ \sigma_{ik})^{-1}Q 
\otimes Q_{ikl}                                         \\ 
& = & (\sigma_{kl},\sigma_{kl})^{-1}Q\otimes Q_{ijk}\otimes 
Q_{ikl}                                                  \\ 
& = & Q_{ijk}\otimes Q_{ikl}. 
\end{eqnarray*} 
Similarly we have $((\sigma_{kl}\circ \sigma_{jk})
\circ \sigma_{ij},\sigma_{il})^{-1}Q = Q_{jkl}\otimes 
Q_{ijl}$.  Also it is clear that 
$$
((\sigma_{kl}\circ \sigma_{jk})\circ \sigma_{ij},
\sigma_{kl}\circ (\sigma_{jk}\circ \sigma_{ij}))^{-1}Q 
= A_{ijkl}, 
$$
where we denote the pullback bundle $(s_i,s_j,s_k,
s_l)^{-1}A$ on $U_{ijkl}$ by $A_{ijkl}$.  It follows as 
above that there is an isomorphism 
$$
A_{ijkl} = Q_{jkl}\otimes Q_{ikl}^*\otimes Q_{ijl}
\otimes Q_{ijk}^* 
$$
of $\cstar$ bundles on $U_{ijkl}$.  Choose a section 
$\rho_{ijk}$ of $Q_{ijk}$ over $U_{ijk}$ and define 
a map $\epsilon_{ijkl}\colon U_{ijkl}\to \cstar$ by 
$\rho_{jkl}\otimes \rho_{ikl}^*\otimes \rho_{ijl}
\otimes \rho_{ijk}^* = a(s_i,s_j,s_k,s_l)\epsilon_
{ijkl}$.  As above $\epsilon_{ijkl}$ satisfies 
the \v{C}ech $3$-cocycle condition $\d(\epsilon)_{ijklm} 
= 1$ and hence is a representative of a class in 
$\check{H}^3 (M;\underline{\C}^{\times}_M) = H^4(M;\Z)$.  
It is straightforward to check that these two methods 
of assigning a \v{C}ech $3$-cocycle to a bundle $2$-gerbe 
give rise to the same class in $H^4(M;\Z)$.  

It is also a straightforward exercise to define such 
notions as the pullback of a bundle $2$-gerbe and the 
product of two bundle $2$-gerbes and prove that the 
four classes behave as one would expect under these operations. 
   
\section{Bundle $2$-gerbe Connections and $2$-curvings} 
\label{sec:eight} 

Just as there is a notion of a bundle gerbe connection 
on a bundle gerbe, there is also a notion of a bundle 
$2$-gerbe connection on a bundle $2$-gerbe $(Q,Y,X,M)$.  
This requires a choice of both a bundle gerbe connection 
$\nabla$ on the bundle gerbe $(Q,Y,X^{[2]})$ and a 
curving $f$ for $\nabla$.  

\begin{definition}[\cite{Ste}]  
\label{def:bundle $2$-gerbe connection} 
Let $(Q,Y,X,M)$ be a bundle $2$-gerbe.  A \emph{bundle 
$2$-gerbe connection} on $Q$ is a pair $(\nabla,f_1)$ 
where $\nabla$ is a bundle gerbe connection on the 
bundle gerbe $Q$ and $f_1$ is a curving for $\nabla$ such 
that the associated $3$-curvature $\omega$ on $X^{[2]}$ 
satisfies $\d(\omega) = 0$.  
\end{definition} 
  
For a proof that bundle gerbe connections always 
exist, see \cite{Ste}.  Note that this is a non-trivial 
fact to prove, as one has to deal with two complexes 
$(\Omega^{p}(X^{[\bullet]}),\d_X)$ and 
$(\Omega^{p}(Y^{[\bullet]}),\d_Y)$ associated to the 
two local-section-admitting surjections $\pi_X\colon 
X\to M$ and $\pi_Y\colon Y\to X^{[2]}$.  The idea of the 
proof is to first choose any bundle gerbe connection 
$\nabla$ on $Q$ and any curving $f$ for $\nabla$.  
Then one can show that there is a two form $\mu\in 
\Omega^{2}(X^{[3]})$ such that $\d(\omega) = d\mu$, 
where $\omega$ is the $3$-curvature associated to 
the bundle gerbe connection $\nabla$ and curving $f$.  
Similarly one can show that there is a one form 
$\a\in\Omega^1(X^{[4]})$ such that $\d(\mu) = d\a$ and 
moreover $\d(\a) = 0$.  Hence, using the exactness 
of the complex~(\ref{eq:exact complex of forms}) one can solve the 
equation $\a = \d(\b)$ for some one form $\b \in 
\Omega^1(X^{[3]})$.  Continuing in this way one can show 
that it is possible to adjust the curving $f$ by the pullback 
of a two form on $X^{[2]}$ so that 
the $3$-curvature $\omega'$ associated 
to $\nabla$ and the new curving $f_1$ satisfies $\d(\omega') = 0$.  

Given a bundle $2$-gerbe connection $(\nabla,f_1)$ on a 
bundle $2$-gerbe $(Q,Y,X,M)$ we can solve the equation 
$\omega = \d(f_2)$ for some three form $f_2$.  A choice 
of $f_2$ is called a \emph{$2$-curving} for the bundle 
$2$-gerbe connection $(\nabla,f_1)$.  Given a choice 
of $2$-curving $f_2$, we have $\d(df_2) = 0$ and hence 
$df_2 = \pi^*(\Theta)$ for some necessarily closed 
four form $\Theta$ on $M$.  We call $\Theta$ the \emph{
four curvature} of the bundle $2$-gerbe connection $(\nabla,f_1)$ 
and $2$-curving $f_2$.  We have the following Proposition.  

\begin{proposition}[\cite{Ste}]  
\label{prop:four curvature} 
The four curvature $\Theta$ is a closed, integral 
four form on $M$ which represents the image in 
$H^4(M;\Reals)$ of the class in $H^4(M;\Z)$ represented 
by the \v{C}ech cocycle $g_{ijkl}$.   
\end{proposition} 
  
As an example of this structure, consider the tautological 
bundle $2$-gerbe on a $3$-connected manifold $M$ associated 
to a closed, integral four form $\Theta$ on $M$.  Recall 
that the tautological bundle $2$-gerbe $(Q,Y,\mathcal{P}M,M)$ 
was defined by constructing the tautological bundle gerbe 
on each fiber of $\pi\colon \mathcal{P}M\to M$.  Another way 
of viewing this construction is to first pull back the four 
form $\Theta$ on $M$ to $\mathcal{P}M$.  Since $\mathcal{P}M$ 
is contractible we can solve $\pi^* \Theta = df_2$ for some 
three form $f_2$ on $\mathcal{P}M$.  Then it is easy to 
see that the three form $\d(f_2)$ on $\mathcal{P}M^{[2]}$ 
is closed.  Since $M$ is $3$-connected, $\mathcal{P}M^{[2]}$ 
is $2$-connected and we can construct the tautological bundle gerbe 
$(Q,Y,\mathcal{P}M^{[2]})$ on $\mathcal{P}M^{[2]}$ 
from the three form $\d(f_2)$ 
using the methods of \cite{Mur} and \cite{CarMurWan}.  
$(Q,Y,\mathcal{P}M,M)$ is then the tautological bundle $2$-gerbe.  In 
\cite{Mur} it is shown how to construct a bundle gerbe connection 
on the tautological bundle gerbe over $\mathcal{P}M^{[2]}$ 
and a curving such that the associated $3$-curvature is $\d(f_2)$.  
This choice of bundle gerbe connection and curving therefore 
defines a bundle $2$-gerbe connection on the tautological bundle 
$2$-gerbe and $f_2$ provides a $2$-curving for this bundle $2$-gerbe 
connection.  $\Theta$ is then the associated $4$-curvature.       

It can be shown \cite{Ste} that given a bundle $2$-gerbe 
$(Q,Y,X,M)$ with bundle $2$-gerbe connection $(\nabla,
f_1)$ and $2$-curving $f_2$ there is a class $D(Q,\nabla,
f_1,f_2)$ in the Deligne hypercohomology group 
$H^3(M;\underline{\C}^{\times}_M\to \underline{\Omega}^1_M 
\to \underline{\Omega}^2_M\to \underline{\Omega}^3_M)$ 
associated to $Q$.  As a consequence of this one 
can show that the class in $H^4(M)$ 
defined by the $4$-curvature $\Theta$ equals the 
image in $H^4(M)$ of the class in $H^3(M;\underline{\C}^\times_M)$ 
defined by the \v{C}ech 
$3$-cocycle $g_{ijkl}$.   

\section{Bundle 2-Gerbes and the First Pontraygin Class} 
\label{sec:nine} 

Suppose we are given a principal $G$ 
bundle $P \to  M$, where $G$ is a 
compact, simply connected, simple Lie group.  
Then it is well known that  
$\pi_{2}(G) = 0$ and $H^{3}(G;\Z) = \Z$.  
It is shown in \cite{Bry} that there 
is a closed, bi-invariant three form $\nu$ on 
$G$ with integral periods which 
represents the canonical generator of 
$H^{3}(G;\Z)=\Z$.  If $G = SU(N)$, then 
$\nu$ is the three form $\frac{1}{24\pi^{2}}
\text{tr}(dgg^{-1})^{3}$.    

Recall from \cite{CarMurWan} that 
we can define a bundle gerbe 
$(Q,\P G,G)$ on $G$ with three 
curvature equal to $\nu$.  
The fibre 
of $Q\to \P G^{[2]}$ at a 
point $(\a,\b) \in \P G^{[2]}$ is the set of 
all equivalence classes $[\phi,z]$ where $z\in \cstar$ and   
$\phi\colon I^2 \to G$ is a homotopy with end 
points fixed between $\a$ and $\b$.  Two pairs 
$(\phi_1,z_1)$ and $(\phi_2,z_2)$ are declared 
equivalent if for all homotopies $F\colon I^3 
\to G$ with end points fixed between $\phi_1$ and 
$\phi_2$ we have $z_2 = z_1 \exp(\int_{I^3}F^* \nu)$.   
The bundle gerbe product is defined by 
$$
[\phi_{1},z_{1}]\otimes [\phi_{2},z_{2}] \mapsto 
[\phi_{1}\phi_{2},z_{1}z_{2}],
$$   
where $\phi_{1}\phi_{2}$ denotes 
the homotopy defined by 
$$
(\phi_{1}\phi_{2})(s,t) = \begin{cases} 
                        \phi_{1}(2s,t) & \text{for $0 \leq s \leq 1/2$} \\
                        \phi_{2}(2s-1,t) &\text{for $1/2 \leq s \leq 1$} \\
                        \end{cases}     
$$
It is shown in \cite{CarMurWan} that 
this is well defined, associative, etc.  

\begin{proposition}[\cite{Ste}]  
The bundle gerbe $(Q,\P G,G)$ is a 
simplicial bundle gerbe on the 
simplicial manifold $NG$.
\end{proposition} 

\begin{proof} 
We first need to define the bundle gerbe 
morphism $m = (\hat{m},m,\text{id})$ 
which maps 
$$
m:d_{0}^{-1}Q\otimes d_{2}^{-1}Q 
\to d_{1}^{-1}Q.
$$
Define 
$m:d_{0}^{-1}\P G\times_{G^{2}}d_{2}^{-1}\P G\to d_{1}^{-1}\P G$ 
covering the identity on $G^{2} = G\times G$ 
by sending $(\a,\b)$ to the piecewise smooth path 
$\a\circ \a(1)\b$ given by 
$$
(\a\circ \a(1)\b)(t) = \begin{cases} 
                    \a(2t), & 0 \leq t \leq 1/2   \\ 
                    \a(1)\b(2t-1), & 1/2\leq t \leq 1. \\ 
                    \end{cases}  
$$
Next, we need to define a $\cstar$ 
equivariant map 
$\hat{m}:d_{0}^{-1}Q\otimes d_{2}^{-1}Q\to d_{1}^{-1}Q$ 
covering 
$$
m^{[2]}:(d_{0}^{-1}\P G\times_{G^{2}}
d_{1}^{-1}\P G)^{[2]}\to d_{1}^{-1}\P G^{[2]}
$$ 
and check that it commutes with the 
bundle gerbe product.  
So take pairs $(\phi, z)$ and 
$(\psi, w)$ where $z,w \in \cstar$ 
and $\phi:I^{2}\to G$ and 
$\psi:I^{2}\to G$ are homotopies with 
endpoints fixed between paths 
$\a_{1},\a_{2}$ and $\b_{1},\b_{2}$ 
respectively.  Then we put 
$$
\hat{m}((\phi,z),(\psi,w)) = 
(\phi\circ \phi(0,1)\psi ,zw)    
$$
where $\phi\circ \phi(0,1)\psi :I^{2}\to G$ 
is the homotopy with endpoints fixed 
between $\a_{1}\circ \a_{1}(1)\b_{1}$ and 
$\a_{2}\circ \a_{2}(1)\b_{2}$ given by 
$$
(\phi\circ \phi(0,1)\psi)(s,t) = 
                     \begin{cases}  
                      \phi(s,2t), & 0\leq t \leq 1/2    \\ 
                      \phi(0,1)\psi(s,2t-1), & 1/2 \leq t \leq 1. \\ 
                      \end{cases} 
$$
We need to check firstly that this 
map is well defined --- that is it respects 
the equivalence relation 
$\equivalencerelation$ --- and secondly that 
$\hat{m}$ commutes with the bundle 
gerbe products.  
So suppose 
$(\phi,z)\equivalencerelation (\phi^{'},z^{'})$ and 
$(\psi,w)\equivalencerelation (\psi^{'},w^{'})$, 
where $\phi$ and $\phi^{'}$ are homotopies 
with endpoints fixed between paths $\a_{1}$ and 
$\a_{2}$ and where $\psi$ and $\psi^{'}$ are 
homotopies with endpoints fixed between 
paths $\b_{1}$ and $\b_{2}$.  We want to show 
that 
$$  
(\phi\circ \phi(0,1)\psi,zw)
\equivalencerelation 
(\phi^{'}\circ \phi^{'}(0,1)\psi^{'},z^{'}w^{'}).   
$$  
Therefore we want to show that for all 
homotopies $H:I^{3}\to G$ with endpoints fixed 
between $\phi\circ \phi(0,1)\psi$ and 
$\phi^{'}\circ \phi^{'}(0,1)\psi^{'}$ we have 
$$
z^{'}w^{'}
 = zw
\exp (\int_{I^{3}}H^{*}\nu).
$$
Note that if $\Phi:I^{3}\to G$ is a 
homotopy with endpoints fixed between 
$\phi$ and $\phi^{'}$ and  
$\Psi:I^{3}\to G$ is a homotopy with 
endpoints fixed between $\psi$ and 
$\psi^{'}$, then by integrality 
of $\nu$ we have  
$$
\exp (\int_{I^{3}}H^{*}\nu) = 
\exp (\int_{I^{3}}(\Phi\circ \Phi(0,0,1)\Psi)^{*}\nu). 
$$
Therefore we are reduced to showing 
that 
$$  
z^{'}w^{'} 
= zw \exp (\int_{I^{3}}(\Phi\circ \Phi(0,0,1)\Psi)^{*}\nu).   
$$ 
We have 
$$
\exp (\int_{I^{3}}(\Phi\circ \Phi(0,0,1)\Psi)^{*}\nu) = 
\exp (\int_{I^{3}}\Phi^{*}\nu) \exp (\int_{I^{3}}
(\Phi(0,0,1)\Psi)^{*}\nu). 
$$
By the bi-invariantness of $\nu$, we get 
$(\Phi(0,0,1)\Psi)^{*}\nu = \Psi^{*}\nu$, hence 
$$
\exp (\int_{I^{3}}(\Phi * \Phi(0,0,1)\Psi)^{*}\nu ) 
= \exp (\int_{I^{3}}\Phi^{*}\nu)\exp (\int_{I^{3}}\Psi^{*}\nu), 
$$
which implies the result.  
Hence $\hat{m}$ is well defined.  It is 
a straightforward matter to verify that 
$\hat{m}$ respects the bundle gerbe products.  

It remains to show that there is a 
transformation of the bundle gerbe 
morphisms $m_{1}$ and 
$m_{2}$ over $G\times G\times G$ 
which satisfies the compatibility 
criterion over $G\times G\times G\times G$.  
This has already been done above for 
the tautological bundle $2$-gerbe and the 
proof given there carries over to this case.  
\end{proof} 

Suppose that we have a principal 
$G$ bundle $\pi:P\to M$.   
Form the canonical map $\tau\colon 
P^{[2]}\to G$ defined by $p_2 = 
p_1\tau(p_1,p_2)$ for $p_1$ and $p_2$ 
in the same fiber.  We can extend $\tau$ to 
define maps $\tau\colon P^{[q]}\to G^{q-1}$ for 
any $q\geq 2$ by 
$$
\tau(p_1,p_2,\ldots,p_q) = (\tau(p_1,p_2),\ldots, 
\tau(p_{q-1},p_q)).  
$$
Notice that $\tau$ defines a simplicial map 
$P^{[\bullet]} \to NG_{\bullet}$ between the simplicial 
manifolds $P^{[\bullet]}$ and $NG_{\bullet}$.  
Clearly the pullback bundle gerbe $\tau^{-1}Q = 
(\tilde{Q},\tilde{P},P^{[2]})$ is a 
simplicial bundle gerbe on the simplicial 
manifold $P^{[\bullet]}$.  We have the following 
Proposition.  
 
\begin{proposition} 
\label{prop:Chern-Simons bundle 2-gerbe} 
The quadruple of manifolds $(\tilde{Q},\tilde{P},P,M)$ 
is a bundle $2$-gerbe.  
\end{proposition} 

Brylinski and McLaughlin \cite{BryMac} \cite{BryMac1} 
defined a canonical $2$-gerbe associated to a 
principal $G$ bundle $P$ on $M$ where $G$ was a 
compact, simple, simply connected Lie group.  They 
showed that the class in $H^4(M;\Z)$ associated to 
this $2$-gerbe was equal to the first Pontryagin 
class $p_1$ of the principal bundle $P$.  If we 
calculate the four class associated to the bundle 
$2$-gerbe $\tilde{Q}$ of 
Proposition~\ref{prop:Chern-Simons bundle 2-gerbe} 
then we recapture the result of Brylinski and McLaughlin.  

\begin{proposition}[\cite{BryMac},\cite{BryMac1}]  
\label{prop:1st Pontraygin class}    
The class in $H^{4}(M;\Z)$ associated 
to the bundle $2$-gerbe $\tilde{Q}$ is the transgression of 
$[\nu]$, that is the first Pontryagin 
class $p_{1}$ of $P$.  
\end{proposition} 

\begin{proof} 
We will calculate the \v{C}ech four class of the 
bundle $2$-gerbe $\tilde{Q}$ and 
show that it is exactly equal to the \v{C}ech cocycle 
obtained by Brylinski and McLaughlin in 
\cite{BryMac} and \cite{BryMac1}.  We then 
apply Theorem 6.2 of \cite{BryMac} to conclude 
that this \v{C}ech four class is $p_{1}$.  We  
calculate the \v{C}ech cocycle $g_{ijkl}$ as follows.  
First choose an open cover $\{U_{i}\}_{i\in I}$ 
of $M$ relative to which $\pi\colon P\to M$ has 
local sections $s_{i}$.  Since $\tilde{P}\to 
P^{[2]}$ is a fibration, we can choose sections $\sigma_{ij}$  
of the pullback fibration $(s_i,s_j)^{-1}\tilde{P}  
= \tilde{P}_{ij}\to U_{ij}$.   
This is 
equivalent to choosing maps $\c_{ij}\colon U_{ij}
\times I\to G$ such that $\c_{ij}(m,0) = 
1$ and $\c_{ij}(m,1) = g_{ij}(m)$.  Next 
we choose sections $\rho_{ijk}\colon U_{ijk}
\to (\sigma_{ik},\sigma_{jk}\circ \sigma_{ij})
^{-1}\tilde{Q}_{ik}$.  This amounts to 
choosing maps $\c_{ijk}\colon U_{ijk}\times I
\times I\to G$ such that $\c_{ijk}(m,0,t) = 
\c_{ik}(m,t)$, $\c_{ijk}(m,1,t) = 
(\c_{ij}\circ g_{ij}\c_{jk})(m,t)$, 
$\c_{ijk}(m,s,0) = 1$ and $\c_{ijk}(m,s,1) 
= g_{ij}(m)g_{jk}(m)$.  Such maps $\c_{ijk}$ 
exist because $G$ is simply connected.  
Define a section $t_{ijkl}$ of the 
bundle $(\sigma_{il},\sigma_{kl}\circ 
(\sigma_{jk}\circ \sigma_{ij}))^{-1}\tilde{Q}
_{il}$ by $t_{ijkl} = 
(e(\sigma_{kl})\circ \rho_{ijk})
\rho_{ikl}$.   
In a similar manner construct a 
section $s_{ijkl} =  
(\rho_{jkl}\circ e(\sigma_{ij})) 
\rho_{ijl}$ of the bundle $(\sigma_{il},
\sigma_{kl}\circ (\sigma_{jk}\circ \sigma_{ij}))
^{-1}\tilde{Q}_{il}$.  We make $t_{ijkl}$ 
into a section of $(\sigma_{il},(\sigma_{kl}
\circ \sigma_{jk})\circ \sigma_{ij})^{-1}
\tilde{Q}_{il}$ by using the associator 
section: $a(\sigma_{kl},
\sigma_{jk},\sigma_{ij}) t_{ijkl}$.    
Finally we define the cocycle $g_{ijkl}$ 
by $s_{ijkl} = a(\sigma_{kl},
\sigma_{jk},\sigma_{ij}) t_{ijkl}
\cdot g_{ijkl}$.  We can get an explicit formula 
for $g_{ijkl}$ as follows: we choose a homotopy 
with endpoints fixed $H_{ijkl}\colon U_{ijkl}\times 
I\times I\times I\to G$ such that $H_{ijkl}(m,0,s,t) = 
(\c_{ij}\circ g_{ij}\c_{jkl})\c_{ijl}(s,t)$, $H_{ijkl}
(m,1,s,t) =  \bar{a}(\c_{kl},\c_{jk},\c_{ij})
(\c_{ijk}\circ g_{ik}\c_{kl})\c_{ikl}(s,t)$, 
$H_{ijkl}(m,r,0,t) = \c_{il}(m,t)$, $H_{ijkl}
(m,r,1,t) = (\c_{ij}\circ g_{jk}\c_{jk})
\circ g_{ik}\c_{kl}$, $H_{ijkl}(m,r,s,0) = 1$  
and $H_{ijkl}(m,r,s,1) = g_{il}$ and we set 
$g_{ijkl} = \exp(\int_{I^{3}}H_{ijkl}^{*}\nu)$.  
This is just the integral of $\nu$ over the 
tetrahedron shown in the following diagram,  
$$ 
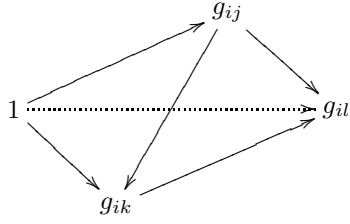
\diagram 
& & g_{ij} \ar[ddl] \ar[dr] &      \\ 
1 \ar[urr] \ar[dr] \ar @1{.>}[rrr] &
& & g_{il}                            \\ 
& g_{ik} \ar[urr] &           &       \\ 
\enddiagram 
$$
as described in \cite{BryMac} and 
\cite{BryMac1}.  Thus our cocycle agrees  
with the cocycle defined by Brylinski and 
McLaughlin.    
\end{proof}  

\section{Higher Gluing Laws} 
\label{sec:ten} 

As a prelude to the discussion of 
trivial bundle $2$-gerbes in the next section, 
we will discuss some features of 
\emph{$2$-descent} (see \cite{Bre}).  We have 
already seen that if we are given a family 
of $\cstar$ bundles $P_i$ defined on an open cover 
$\{U_i\}_{i\in I}$ of a manifold $M$ such that 
there exist isomorphisms $\phi_{ij}\colon 
P_i\to P_j$ satisfying the cocycle condition 
$\phi_{jk}\circ \phi_{ij} = \phi_{ik}$ then 
we can construct a $\cstar$ bundle $P$ 
defined on $M$ which is locally isomorphic 
to $P_i$ over each open set $U_i$.  If we  
replace $\cstar$ bundles by bundle gerbes 
then new complications arise.  Rather than demanding that 
the equation $\phi_{jk}\circ \phi_{ij} = \phi_{ik}$ 
is satisfied on the nose, we can settle for the weaker 
condition that there is a transformation of bundle gerbe 
morphisms $\psi_{ijk}\colon \phi_{jk}\circ \phi_{ij}
\Rightarrow \phi_{ik}$ which satisfies a certain 
cocycle condition --- the \emph{non-abelian $2$-cocycle 
condition}.  We shall see that it is still possible 
to `glue' the various bundle gerbes $P_i$ together 
to form a bundle gerbe $P$ on $M$.    
 
Suppose we are given an open cover 
$\{U_i\}_{i\in I}$ of a manifold 
$M$ such that there exist bundle 
gerbes $Q_i$ over $U_i$.  Suppose also 
that over each intersection $U_{ij}$ there 
exist bundle gerbe morphisms $\phi_{ij}
\colon Q_i|_{U_{ij}}\to Q_j|_{U_{ij}}$.   
Suppose as well, that over each triple 
intersection $U_{ijk}$ 
there exist transformations of bundle gerbe 
morphisms $\psi_{ijk}\colon \phi_{jk}|_{U_{ijk}}
\circ \phi_{ij}|_{U_{ij}}\Rightarrow 
\phi_{ik}|_{U_{ijk}}$.  Finally, suppose that 
the diagram of transformations 
of bundle gerbe morphisms in Figure~\ref{fig:non-abelian 
2-cocycle condition} commutes. 
\begin{figure} 
$$
\xymatrix{ 
Q_i \ar[rr]^-{\phi_{ij}} \ar[dd]_-{\phi_{il}} 
\ar[ddrr] & & Q_j \ar[dd]_-{\phi_{jk}} 
\ar @2{->}[dl]_-{\psi_{ijk}} & Q_i \ar[dd]_<<<<<<<{\phi_{il}} 
\ar[rr]^-{\phi_{ij}} & & Q_j \ar[ddll] 
\ar[dd]^-{\phi_{jk}}                                       \\ 
& \ar @2{->}[l]^-{\psi_{ikl}} & \ar @2{-}[r] & & 
\ar @2{->}[l]_-{\psi_{ijl}} &                                   \\ 
Q_l & & Q_k \ar[ll]_-{\phi_{kl}} & Q_l & & Q_k 
\ar @2{->}[ul]^-{\psi_{jkl}} \ar[ll]^-{\phi_{kl}}             } 
$$
\caption{The non-abelian $2$-cocycle condition} 
\label{fig:non-abelian 2-cocycle condition} 
\end{figure}
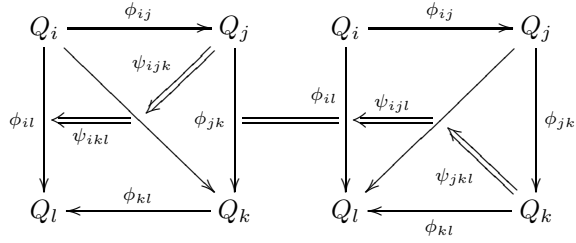 
If we let $L_{ijk}$ denote the $\cstar$ bundle 
$D_{\phi_{jk}\circ \phi_{ij},\phi_{ik}}$ 
on $U_{ijk}$ then we have an isomorphism $L_{ijk}
\otimes L_{ikl} = L_{ijl}\otimes L_{jkl}$ of 
$\cstar$ bundles on $U_{ijkl}$ or, put another way, 
a canonical trivialisation $L_{jkl}\otimes 
L_{ikl}^*\otimes L_{ijl}\otimes L_{ijk}^* = 1$.  The 
condition that the diagram of bundle gerbe transformations 
in Figure~\ref{fig:non-abelian 2-cocycle condition} 
commutes translates into the requirement that the 
induced section $\psi_{jkl}\otimes \psi_{ikl}^*\otimes 
\psi_{ijl}\otimes \psi_{ijk}^*$ matches this canonical 
trivialisation.  We have the following proposition.  

\begin{proposition} 
\label{prop:bundle gerbe gluing} 
Suppose we are given an open cover $\{U_i\}$ of 
$M$ and a triple $(Q_i,\phi_{ij},\psi_{ijk})$ as 
described above.  Then there is a bundle gerbe $Q$ on 
$M$ and bundle gerbe morphisms $\chi_i\colon Q|_{U_i}
\to Q_i$ over $U_i$ together with transformations 
$\xi_{ij}\colon \phi_{ij}\circ \chi_i\Rightarrow 
\chi_j$ which are compatible with the transformations 
$\psi_{ijk}$.       
\end{proposition} 

\begin{proof} 
Suppose the bundle gerbes $Q_i$ are given by triples 
$(Q_i,X_i,U_i)$.  We first construct the bundle gerbe 
$(Q,X,M)$.  Let $X = \coprod_{i\in I}X_i$.  Then the 
fibre product of $X$ with itself over $M$ is 
$X^{[2]} = \coprod_{i,j\in I}X_i\times_M 
X_j$.  Suppose the bundle gerbe morphisms $\phi_{ij}$ 
are given by $\phi_{ij} = (\hat{\phi}_{ij},
\phi_{ij})$.  Define a map $f_{ij}\colon X_i\times_M
X_j\to X_j^{[2]}$ by sending $(x_i,x_j)\in X_i\times_M
X_j$ to $(\phi_{ij}(x_i),x_j)$.  Let $Q_{ij} = 
f_{ij}^{-1}Q_j$.  Define a $\cstar$ bundle $Q$ on 
$X^{[2]}$ by setting $Q = \coprod_{i,j\in I}Q_{ij}$ 
with projection map $Q\to X^{[2]}$ induced by the 
various projections $Q_{ij}\to X_i\times_M X_j$.  
We want to show that the triple $(Q,X,M)$ is a 
bundle gerbe.  We first define the product in $Q$.  
This is a $\cstar$ bundle isomorphism $\pi_1^{-1}Q
\otimes \pi_3^{-1}Q\to \pi_2^{-1}Q$ covering the 
identity on $X^{[3]}$ which satisfies an associativity 
condition on $X^{[4]}$.  Since $X^{[3]} = X_i\times_M
X_j\times_M X_k$ this amounts to finding a $\cstar$ 
bundle map $Q_{jk}\otimes Q_{ij}\to Q_{ik}$ satisfying an 
associativity condition over $X_i\times_M X_j\times_M 
X_k\times_M X_l$.   

Let $u_{jk}\in (Q_{jk})_{(x_j,x_k)}$, $u_{ij} \in 
(Q_{ij})_{(x_i,x_j)}$ for $(x_i,x_j,x_k)\in 
X_i\times_M X_j\times_M X_k$.  Then $u_{jk}\in 
(Q_{j})_{(\phi_{ij}(x_i),x_j)}$ and $u_{ij}\in 
(Q_{k})_{(\phi_{jk}(x_j),x_k)}$.  Apply $\hat{\phi}
_{ij}$ to $u_{ij}$.  Then $\hat{\phi}_{jk}(u_{ij}) 
\in (Q_{k})_{(\phi_{jk}(\phi_{ij}(x_i)),\phi_{jk}
(x_j))}$.  Using the bundle gerbe product in $Q_k$ 
we have that 
$$
u_{jk}\hat{\phi}_{jk}(u_{ij}) \in (Q_{k})_{(\phi_{jk}
(\phi_{ij}(x_i)),x_k)}.  
$$
Let $\hat{\psi}_{ijk}$ denote the section of the $\cstar$ 
bundle $(\phi_{jk}\circ \phi_{ij},\phi_{ik})^{-1}
Q_k$ on $X_i|_{U_{ijk}}$ which descends to $\psi_{ijk}$.  
Using the bundle gerbe product in $Q_k$ again, we have 
that 
$$
u_{jk}\hat{\phi}_{jk}(u_{ij})\hat{\psi}_{ijk}^{-1}(x_i)  
\in (Q_{k})_{(\phi_{ik}(x_i),x_k)}.   
$$
We define a product in $Q$ by sending $u_{jk}\otimes 
u_{ij}$ to $u_{jk}\cdot u_{ij} = u_{jk}\hat{\phi}_{ij}
(u_{ij})\hat{\psi}_{ijk}^{-1}(x_i)$.   
We have to check that this product is 
associative.  This follows easily from 
the following equation satisfied by  
$\hat{\psi}_{ijk}$:  
$$ 
\hat{\psi}_{ikl}(x_i)\hat{\phi}_{kl}
(\hat{\psi}_{ijk}(x_i)) = \hat{\psi}_{ijl}(x_i)
\hat{\psi}_{jkl}(\phi_{ij}(x_i)). 
$$
This equation is a consequence of the coherency 
condition satisfied by $\psi_{ijk}$.  Therefore 
$(Q,X,M)$ is a bundle gerbe.  We now need to define 
the bundle gerbe morphism $Q|_{U_i}\to Q_i$.  
First of all we define a map $X|_{U_i}\to X_i$ 
covering the identity on $U_i$.  If $x_j\in X_j$ 
and $\pi_{X_j}(x_j)\in U_i$, then $\phi_{ji}
(x_j)\in X_i$.  Since $X|_{U_i} = \coprod_{j \in J}
X_j|_{U_i}$ this defines a map $X|_{U_i}\to X_i$.  
Now suppose $(x_j,x_{j'})\in X|_{U_i}^{[2]}$ and 
$u_{jj'}\in Q_{(x_j,x_{j'})}$.  So $u_{jj'}\in 
(Q_{j'})_{(\phi_{jj'}(x_j),x_{j'})}$.  Hence applying 
$\hat{\phi}_{j'i}$ to $u_{ij}$ means that 
$\hat{\phi}_{j'i}(u_{jj'}) \in (Q_{i})_{(\phi_{j'i}(
\phi_{jj'}(x_j)),\phi_{j'i}(x_{j'}))}$.  
Therefore 
$$
\hat{\phi}_{j'i}(u_{jj'})\hat{\psi}^{-1}_{jj'i}(x_j) 
\in (Q_{i})_{(\phi_{ij}(x_j),\phi_{ij'}(x_j'))}.  
$$
This defines a $\cstar$ bundle map $Q|_{U_i}\to Q_i$.  
It is not hard to check that this map commutes with 
the bundle gerbe products on $Q$ and $Q_i$ and 
hence defines a bundle gerbe morphism 
$\chi_i\colon Q|_{U_i}\to Q_i$.  Similarly, 
one can define a transformation of bundle gerbe 
morphisms $\xi_{ij}\colon \phi_{ij}\circ \chi_i \Rightarrow 
\chi_j$ which is compatible with $\psi_{ijk}$.   
\end{proof} 

The triple $(Q_i,\phi_{ij},\psi_{ijk})$ is 
called \emph{$2$-descent data} relative to the 
open covering $\U = \{U_i\}_{i\in I}$.  
One can think of these 
$2$-descent data as being objects of a $2$-category 
$\textbf{2-Desc}(\U)$.  Let  
$(Q_i,\phi_{ij},\psi_{ijk})$ and 
$(P_i,\tilde{\phi}_{ij},\tilde{\psi}_{ijk})$ 
be two sets of $2$-descent data.  A $1$-arrow from 
$(Q_i,\phi_{ij},\psi_{ijk})$ to $(P_i,
\tilde{\phi}_{ij},\tilde{\psi}_{ijk})$ is a pair 
$(f_{ij},\tau_{ij})$ where $f_i\colon Q_i\to P_i$ is a 
bundle gerbe morphism and $\tau_{ij}$ is a transformation 
of bundle gerbe morphisms as pictured in the 
following diagram 
$$
\xymatrix{ 
Q_i \ar[d]_-{\phi_{ij}} \ar[r]^-{f_i} & P_i 
\ar[d]^-{\tilde{\phi}_{ij}} \ar @2{->}[dl]^-{\tau_{ij}} \\ 
Q_j \ar[r]_-{f_j} & P_j                                     } 
$$
which is compatible with $\psi_{ijk}$ and 
$\tilde{\psi}_{ijk}$.  Given two $1$-arrows 
$(f_i,\tau_{ij})$ and $(g_i,\rho_{ij})$ a $2$-arrow 
$(f_i,\tau_{ij})\Rightarrow (g_i,\rho_{ij})$ is a 
transformation of bundle gerbe morphisms $\lambda_i
\colon f_i\Rightarrow g_i$ which is compatible 
with $\tau_{ij}$ and $\rho_{ij}$.  Horizontal and vertical 
composition in $\textbf{2-Desc}(\U)$ is defined 
in the obvious manner.    

The gluing procedure of Proposition~
\ref{prop:bundle gerbe gluing} above allows us to define a $2$-functor 
$\textbf{2-Desc}(\U)\to \BGrb_M$.  The 
action of this functor on objects of $\textbf{2-Desc}
(\U)$ is clear: a triple $(Q_i,\phi_{ij},\psi_{ijk})$ 
of $2$-descent data is mapped to the bundle gerbe 
$Q$ of Proposition~\ref{prop:bundle gerbe gluing}.  With 
a little work one can show that a $1$-arrow $(f_i,\tau_{ij})$ 
from $(Q_i,\phi_{ij},\psi_{ijk})$ to $(P_i,\tilde{
\phi}_{ij},\tilde{\psi}_{ijk})$ induces a bundle 
gerbe morphism $f\colon Q\to P$ and that a $2$-arrow 
$\lambda_i\colon (f_i,\tau_{ij})\Rightarrow (g_i,\rho_{ij})$ 
between two $1$-arrows $(f_i,\tau_{ij})$ 
and $(g_i,\rho_{ij})$ induces a transformation 
of bundle gerbe morphisms $\lambda\colon f\Rightarrow g$.  
Both of these constructions are functorial.       
 
Note that bundle gerbe morphisms do not glue together in the 
fashion that one would like.  One would like to say 
that given bundle gerbes $P$ and $Q$ such that 
relative to some open cover $\{U_i\}_{i\in I}$ of 
$M$ there exist local bundle gerbe morphisms $f_i\colon 
P|_{U_i}\to Q|_{U_i}$ together with transformations of 
bundle gerbe morphisms $\tau_{ij}\colon f_i\Rightarrow 
f_j$ which satisfy the cocycle condition $\tau_{jk}
\tau_{ij} = \tau_{ik}$, there exists a bundle gerbe 
morphism $f\colon P\to Q$ locally isomorphic to 
$f_i$.  Unfortunately this is not true; it is 
however true for gerbes.     
        
\section{Trivial Bundle 2-Gerbes} 
\label{sec:eleven} 

In \cite{Mur} it was shown that a bundle 
gerbe $P$ had vanishing Dixmier-Douady class 
precisely when the bundle gerbe was trivial --- 
ie $P$ was of the form $\d(T)$ for some $\cstar$ 
bundle $T$.  We would like to know under what 
conditions the four class of a bundle $2$-gerbe 
is zero.  We will define a certain class of 
bundle $2$-gerbes, \emph{trivial} bundle $2$-gerbes 
and show that the four class associated to 
a bundle $2$-gerbe belonging to this class vanishes.  
We will then prove that the converse is true.  

\begin{definition} 
\label{def:trivial bundle 2-gerbe} 
Let $(Q,Y,X,M)$ be a bundle $2$-gerbe.  We say 
that $Q$ is \emph{trivial} if there exists 
a bundle gerbe $(L,Z,X)$ on $X$ together with a 
bundle gerbe morphism $\eta\colon \pi_1^{-1}
L\otimes Q\to \pi_2^{-1}L$ over $X^{[2]}$ and a 
transformation of bundle gerbe morphisms 
$\theta$ as pictured in the following diagram: 
$$
\xymatrix{ 
\pi_1^{-1}\pi_1^{-1}L\otimes \pi_1^{-1}Q
\otimes \pi_3^{-1}Q \ar[d]_-{\pi_1^{-1}\eta
\otimes 1} \ar[rr]^-{1\otimes m} & & 
\pi_1^{-1}\pi_1^{-1}L\otimes \pi_2^{-1}Q 
\ar @2{-}[d]                                       \\ 
\pi_1^{-1}\pi_2^{-1}L\otimes \pi_3^{-1}Q 
\ar @2{-}[d] & & \pi_2^{-1}\pi_1^{-1}L\otimes 
\pi_2^{-1}Q \ar @2{->}[ll]_-{\theta} 
\ar[d]^-{\pi_2^{-1}\eta}                       \\ 
\pi_3^{-1}\pi_1^{-1}L\otimes \pi_3^{-1}Q 
\ar[dr]_-{\pi_3^{-1}\eta} & & \pi_2^{-1}
\pi_2^{-1}L                                         \\ 
& \pi_3^{-1}\pi_2^{-1}L \ar @2{-}[ur].               } 
$$
Let us agree to call $\eta_1 = \pi_2^{-1}
\eta\circ (1\otimes m)$ and 
$\eta_2 = \pi_3^{-1}\eta\circ (\pi_1^{-1}
\eta\otimes 1)$.  Then $\theta$ is a 
section trivialising the $\cstar$ bundle $B = 
D_{\eta_1,\eta_2}$ on $X^{[3]}$.  Moreover there 
is a canonical isomorphism $\d(B) = \pi_1^{-1}B\otimes 
\pi_2^{-1}B^*\otimes \pi_3^{-1}B\otimes \pi_4^{-1}B^* 
= A$ of $\cstar$ bundles over $X^{[4]}$.  As a final 
condition we demand that the induced section 
$\d(\theta) = \pi_1^{-1}\theta \otimes \pi_2^{-1}\theta^* 
\otimes \pi_3^{-1}\theta \otimes \pi_4^{-1}\theta^*$ 
of $\d(B)$ is mapped to $a$ under this isomorphism.  
\end{definition} 
 
Suppose we are now given a bundle $2$-gerbe $(Q,Y,X,M)$ 
with vanishing four class.  We will make the additional 
assumption that $\pi_Y\colon Y\to X^{[2]}$ is a fibration.  
Let $\{U_i\}_{i\in I}$ be an open covering of $M$ all 
of whose finite intersections $U_{i_0}\cap \cdots \cap 
U_{i_p}$ are empty or contractible and such that there 
exist local sections $s_i\colon U_i\to X$ of the surjection 
$\pi\colon X\to M$ over $U_i$.  Define maps $\hat{s}_i
\colon X_i\to X^{[2]}$ where $X_i = \pi^{-1}(U_i)$ by 
$\hat{s}_i(x) = (x,s_i(\pi(x)))$.  Let $(L_i,Z_i,X_i)$ 
denote the pullback of the bundle gerbe $(Q,Y,X^{[2]})$ 
to $X_i$ via the map $\hat{s}_i$.  Then $Z_i\to X_i$ is 
a fibering with fibre $(Z_{i})_x$ at $x\in X_i$ equal 
to $Y_{(x,s_i(\pi(x)))}$.  One can also 
define maps $(s_i,s_j)\colon U_{ij}\to X^{[2]}$ 
in the usual fashion by sending $m\in U_{ij}$ 
to $(s_i(m),s_j(m))\in X^{[2]}$.  Let $Y_{ij}$ 
denote the pullback of the fibration $Y\to 
X^{[2]}$ via this map.   
Choose sections 
$\sigma_{ij}$ of the pullback 
fibering $Y_{ij}\to U_{ij}$.  Now we 
can define maps $\phi_{ij}\colon Z_i\to Z_j$ 
by sending $y_i\in Z_i$ to $m(\sigma_{ij},
y_i)\in Z_j$.  The $\phi_{ij}$ extend to define 
bundle gerbe morphisms $\phi_{ij} 
= (\hat{\phi}_{ij},\phi_{ij})\colon L_i\to L_j$ with 
$\hat{\phi}_{ij}(u_i) = \hat{m}(e(\sigma_{ij})
\otimes u_i)$ where $e$ denotes the identity 
section of the bundle gerbe $(Q,Y,X^{[2]})$.  

We now wish to define transformations $\psi_{ijk}
\colon \phi_{jk}\circ \phi_{ij}
\Rightarrow \phi_{ik}$ satisfying the 
non-abelian $2$-cocycle condition over $X_{ijkl}$.            
To do this, first note that $\hat{a}(\sigma_{jk},
\sigma_{ij},y_i)\in Q_{(m(m(\sigma_{jk},\sigma_{ij}),
y_i),m(\sigma_{jk},m(\sigma_{ij},y_i)))}$ where 
$\hat{a}$ denotes the lift of the associator 
section $a$ to $Y\circ Y\circ Y$.  Also, as in 
Section~\ref{sec:seven}, let $\rho_{ijk}$ denote 
a section of the pullback bundle $(m(\sigma_{jk}
,\sigma_{ij}),\sigma_{ik})^{-1}Q$ over $U_{ijk}$.  
Then $\hat{m}(\rho_{ijk}\otimes e(y_i))\in 
Q_{(m(m(\sigma_{jk},\sigma_{ij}),y_i),m(\sigma_{ik},
y_i))}$.  Therefore 
$$
\hat{\psi}_{ijk}(y_i) = \hat{m}(\rho_{ijk}\otimes 
e(y_i))\hat{a}(\sigma_{jk},\sigma_{ij},y_i)^{-1} 
\in Q_{(\phi_{jk}(\phi_{ij}(y_i)),\phi_{ik}(y_i))}.  
$$
Since the \v{C}ech 3-cocycle $g_{ijkl}$ representing 
the four class of $Q$ is trivial, one can show that 
it is possible to choose $\rho_{ijk}$ so that the 
sections $\hat{\psi}_{ijk}$ defined above satisfy the 
non-abelian $2$-cocycle condition.  Therefore, using 
Proposition~\ref{prop:bundle gerbe gluing}, one 
can form a bundle gerbe $(L,Z,X)$ on $X$ which is 
locally isomorphic to each $(L_i,Z_i,X_i)$.  

However, more is true.  The bundle gerbes $(L_i,Z_i,
X_i)$ provide local trivialisations of the bundle $2$-gerbe 
$Q$.  To see this, note that the bundle gerbe morphism 
$m\colon \pi_1^{-1}Q\otimes \pi_3^{-1}Q\to 
\pi_2^{-1}Q$ provides a bundle gerbe morphism 
$\eta_i\colon \pi_1^{-1}L_i\otimes Q\to \pi_2^{-1}
L_i$ by sending a point $(y_i,y)$ of 
$\pi_1^{-1}Z_i\times_{X^{[2]}}Y$ to $m(y_i,y) \in 
\pi_2^{-1}Z_i$ and a point $u_i\otimes u$ of 
$\pi_1^{-1}L_i\otimes Q$ to $\hat{m}(u_i\otimes 
u)$.  One can also define transformations 
$\theta_i$ as in Definition~\ref{def:trivial bundle 
2-gerbe} above.  It is possible to show \cite{Ste}, 
although it is very tedious, that $\eta_i$ 
and $\theta_i$ are compatible with the $2$-descent 
data $(L_i,\phi_{ij},\psi_{ijk})$ relative 
to the open covering $\{X_i\}_{i\in I}$ of $X$.  
It follows that $\eta_i$ and $\theta_{i}$ 
glue together to form a bundle gerbe morphism 
$\eta\colon \pi_1^{-1}L\otimes Q\to \pi_2^{-1}L$ 
and a transformation $\theta$ as in Definition~\ref
{def:trivial bundle 2-gerbe}.  Thus the bundle $2$-gerbe 
$Q$ is trivial.  

One can show \cite{Ste} that it is possible to 
remove the restriction that $\pi_Y\colon Y\to 
X^{[2]}$ be a fibration.  

\begin{proposition}[\cite{Ste}] 
The four class of a bundle $2$-gerbe $(Q,Y,X,M)$ 
vanishes if and only if the bundle $2$-gerbe 
$Q$ is trivial. 
\end{proposition}

\end{document}